\documentclass[11pt]{article}

\usepackage{amsmath, amsthm, amsfonts}
\usepackage{listings}
\usepackage{amsfonts}
\usepackage{amssymb}
\usepackage{float}
\usepackage{anysize}
\usepackage{graphicx}
\usepackage{mathrsfs} 
\usepackage{hyperref}
\usepackage{cite}
\usepackage{color}
\usepackage{subcaption}
\usepackage{color}
\usepackage{cancel}
\usepackage{bbm}
\usepackage{tikz}

\usepackage{fullpage} 

\usepackage{amsmath,latexsym,amssymb,amsthm,array,amsfonts}

\usepackage{mathrsfs,dsfont,bbm}

\usepackage[sectionbib]{natbib}

\usepackage{layout}
\usepackage{layout}
\usepackage{eufrak}
\newtheorem{prop}{Proposition}
\newtheorem{lemma}{Lemma}

\newtheorem{theorem}{Theorem}

\newtheorem{remark}{Remark}

\def\real{{\mathord{{\rm I\kern-2.8pt R}}}}        
\def\inte{{\mathord{{\rm I\kern-2.8pt N}}}}

\def\sZZ{{\rm Z\kern-2.8ptem{}Z}}

\def\z{{\mathchoice
  {\sZZ}
  {\sZZ}
  {\rm Z\kern-0.30em{}Z}
  {\rm Z\kern-0.25em{}Z} }}
\def\sQQ{{\kern 0.27em \vrule height1.45ex width0.03em depth0em
          \kern-0.30em \rm Q}}
\def\qu{{\mathchoice
    {\sQQ}
    {\sQQ}
  {\kern 0.225em \vrule height1.05ex width0.025em depth0em \kern-0.25em \rm Q}
  {\kern 0.180em \vrule height0.78ex width0.020em depth0em \kern-0.20em \rm Q}
        }}
\def\sCC{{\kern 0.27em \vrule height1.45ex width0.03em depth0em
          \kern-0.30em \rm C}}
\def\complex{{\mathchoice
    {\sCC}
    {\sCC}
  {\kern 0.225em \vrule height1.05ex width0.025em depth0em \kern-0.25em \rm C}
  {\kern 0.180em \vrule height0.78ex width0.020em depth0em \kern-0.20em \rm C}
        }}

\newcommand{\E}{\mathbb{E}}

\newcommand{\ba}{\begin{array}}
\newcommand{\ea}{\end{array}}
\newcommand{\be}{\begin{equation}}
\newcommand{\ee}{\end{equation}}
\newcommand{\bea}{\begin{eqnarray}}
\newcommand{\eea}{\end{eqnarray}}
\newcommand{\beaa}{\begin{eqnarray*}}
\newcommand{\eeaa}{\end{eqnarray*}}

\newcommand{\eps}{\varepsilon}

\def\b{\beta}

\def\z{\zeta}

\font\tenmath=msbm10 \font\sevenmath=msbm7 \font\fivemath=msbm5
\newfam\mathfam \textfont\mathfam=\tenmath
\scriptfont\mathfam=\sevenmath \scriptscriptfont\mathfam=\fivemath

\def \b{\noindent}

\def \={{\buildrel {\rm (law)} \over =}}

\def\qed{ \hfill \vrule width.25cm height.25cm depth0cm\smallskip}

\newcommand{\basa}{\begin{assumption}}
\newcommand{\easa}{\end{assumption}}

\newcommand{\bas}{\begin{assum}}
\newcommand{\eas}{\end{assum}}

\newcommand{\ignore}[1]{}
\begin{document}

\renewcommand{\thefootnote}{\fnsymbol{footnote}}

\renewcommand{\thefootnote}{\fnsymbol{footnote}}

\title{Limit distribution of the least square estimator with observations sampled at random times driven by standard Brownian motion}
\author{ Tania Roa \thanks{CIMFAV, Facultad de Ingenier\'ia, Universidad de Valpara\'iso. Email:tania.roa@postgrado.uv.cl} \and Soledad Torres\thanks{CIMFAV, Facultad de Ingenier\'ia, Universidad de Valpara\'iso. Email: soledad.torres@uv.cl} \and Ciprian Tudor \thanks{ Laboratoire Paul Painlev\'e,  Universit\'e de Lille 1, F-59655 Villeneuve d'Ascq,  France. Email: ciprian.tudor@univ-lille.fr}
}

\maketitle

\begin{abstract}

 In this article, we study the limit distribution of the least square estimator, properly normalized, from a regression model in which observations are assumed to be finite ($\alpha N$) and sampled under two different random times. Based on the limit behavior of the characteristic function and convergence result we prove the asymptotic normality for the least square estimator. We present simulations results to illustrate our theoretical results.
\end{abstract}

\vskip0.3cm

{\bf 2010 AMS Classification Numbers:} 60G22; 62J86; 62M09.

\vskip0.3cm


\section{Introduction}

The estimation of the parameters of a stochastic process  on the basis of its random sampling (i.e. the process is observed at random times) received a wide attention in the past. Such a problem is well motivated by practical aspects. Indeed, the measuring instruments (classical or modern, such as sattelites)  may introduce random disturbances to the data. For examples, transaction data in finance arrive in irregular time intervals (see e.g. \cite{ER}), and the same happens with  biological  signals in medicine, such as heart rate (see e.g. \cite{BB}). Other examples of appearances of random models observed at unequally, possibly random, times can be found, among many others,  in climatology (see \cite{C1}, \cite{C2}) or computer science.
Although quite natural, the hypothesis of random sampling  for stochastic models leads to more complex estimators and, in general, some particular choices for the random observation times are considered in the literature. For instance, \cite{DY}, the authors studied  a diffusion process observed at independent Poisson times, \cite{J} the situation when the $i$th observation depends on the previous $i-1$ observations is considered, while \cite{Vi1}, \cite{Vi2} the authors used the so-called jittered and renewal sampling.

Our purpose is to analyse the asymptotic properties of the least squares estimator (LSE in the sequel) for a simple regression model driven by a standard Wiener process, i.e. 

\begin{equation}\label{1}
Y_{\tau_{i+1}} = a\tau _{i+1} + W_{\tau_{i+1}} - W_{\tau_{i}}, \hskip0.5cm i=0,..., N-1
\end{equation}
where $\tau_{i}, i=0,.., N$ are random times, independent of $W$, with $\tau_{0}:=0$. We choose to work with the random sampling proposed by \cite{Vi1}, \cite{Vi2} which includes two types of randomness: the jittered sampling (the observation times are $\frac{i}{N}$, $i=1,..., N$   perturbed by a ``small" uniform random  variable) or the renewal sampling (the $i$th observation times is a sum of $i$ independent positive  random variables, so the randomness is somehow progressive).

We construct  a least squares estimator (LSE) for the drift parameter $a$   of the model   (\ref{1}) and then we analyse its asymptotic properties. Our proofs are based on a sharp calculation of the mean square of the estimator and of its conditional distribution given the random times, which is Gaussian. We also use some results given by \cite{araya2019} where the asymptotic behavior of the denominator of the LSE estimator is obtained. 

We organized our paper as follows. In Section 2 we describe the model and we include a discussion about the number of random observations used to define to estimator. In Section 3 we calculate exactly the mean square norm  of the estimator when the number of observations is large enough while in Section 4 we give the asymptotic distribution for the LSE. Many of our theoretical results are illustrated by numerical simulations in Section 5. 

\section{Preliminaries}

Let us now introduced the random times  considered in our model (\ref{1}). Our examples are inspired from \cite{araya2019} and \cite{Vi1}.

\subsection{Random times}

Let  $T=1$ and $ \tau = \lbrace  \tau_{i}; i=0,\ldots,N \rbrace$ a strictly increasing sequence of random points over time, where  $N$ is the last integer such that $\tau_{N-1} \leq 1$, which exhibits one of the following two features.

\begin{enumerate}
\item { \bf Jittered sampling (JS)}. First, we assume that we observed a certain process at regular times  $\tau$ with period $\delta = 1/N >0$ but contaminated by an additive noise  $\nu$ which represents possible measurement errors. Then the sequence of random times $\tau_i, \quad 0 \le i \le N$ is defined as 
\begin{eqnarray} \label{js}
\tau_{i, N}=:\tau_{i} = \dfrac{i}{N} +  \nu_{i,N}, \quad i = 1, \dots, N \; and\; \tau_0 := 0, \dots , 
\end{eqnarray}
where $\lbrace \nu_{i,N}; \quad 1 \le i \le N \rbrace$ constitutes a triangular array of independent and identically distributed set of random variables with common density function depending on $N$, called $g_{N}(t)$, which is assumed to be symmetric in  $\left[ -\frac{1}{2N}, \frac{1}{2N} \right]$ for all $i= 1,\ldots , N$. From now on, we state the following about $\nu_{i,N}$
\begin{itemize}
\item $\mathbb{E} \left[ \nu_{i,N} \right] = 0$, and
\item $\mathbb{E} \left[ \nu^{2}_{i,N} \right]$ satisfies
\begin{equation}\label{2d-1}
\mathbb{E} \left[ \nu^{2}_{i,N} \right]= c_{1} \frac{1}{N ^ {2}} \mbox{ with } c_{1}>0.
\end{equation}

\end{itemize}
Some distributions that satisfy the latter statement are the uniform distribution in $\left[ -\frac{1}{2N} , \frac{1}{2N} \right]$, triangular distribution with parameters $\left( -\frac{1}{2N} , 0 , \frac{1}{2N} \right)$ and the raised cosine distribution with parameters $\mu=0$ and $s=\frac{1}{2N}$. For instance, $c_{1}= \frac{1}{12}$ when $g_{N}$ is the uniform distribution over the interval  $\left[ -\frac{1}{2N} , \frac{1}{2N} \right]$ and $c_{1}=\frac{1}{24} $ when $g_{N}$ is the  triangular distribution with parameters $\left( -\frac{1}{2N} , 0 , \frac{1}{2N} \right)$.

\item { \bf Renewal sampling (RP)}. In this case, the sequence $\tau$ satisfies the following property
\begin{eqnarray}
\tau_{i} = \sum_{j=1}^{i} t_j   \    \    \   \  i=1,2,... \    \    \  \mbox{and} \, \tau_{0} := 0,
\label{rp}
\end{eqnarray}
where $  \lbrace t_j  , 1 \le j \rbrace$  is a sequence of independent and identically distributed random variables, with a common distribution function  $G$ with support in $[0,\infty)$. In this work we consider that $G$ is an  exponential distribution with parameter $\lambda =N$, i.e. it has density function $g(t)= Ne ^ {-Nt} 1_{(0, \infty) }(t).$ The random times $\tau_{i} $  given by   (\ref{rp}) actually depend also on $N$ but we still use the notation $\tau_{i,N}= : \tau_{i}$, for simplicity. 
\end{enumerate}

\subsection{The number of observations}\label{sec22}

Assume that we observe a stochastic process  $Y$ at times $\tau_{1},..., \tau _{ [\alpha N]}$ with $\tau_{i} <\tau _{i+1} $ for every $i\geq 1$ and with  $0<\alpha \leq 1$. We want to ensure that our observation  period  remains, almost surely,  inside the interval $[0, T]$ with $T=1$. That is, we would like to have that the last observation $ \tau _{[\alpha N ]} $ is almost surely less that $1$ for $N$ sufficiently large. In the case of jittered sampling, this is always true for $\alpha =1$, due to our hypothesis (\ref{2d-1}). Indeed, $\tau_{N-1} = \frac{N-1}{N} + \nu _{N-1, N}$ and $\tau_{N} =  1 + \nu _{N, N}$ and then $P (\tau_{N-1} >1)=0$ while $P(\tau_{N} >1)=\frac{1}{2}$, so $\tau_{N-1}$ is  almost surely in the observation interval $[0,1]$. In this case we assume $\alpha=1$ and $Y_{\tau_{N}}=0$.

On the other hand, in the situation of the renewal sampling, we have $\tau _{N} \sim G (N, N) $ (by $(G(a, \lambda)$ we denote the Gamma  law with parameters $a>0, \lambda >0$)) and by a result of \cite{gaut1977}, 
$$\mathbb{P}(\tau_{N}>1)= \mathbb{P}( G(N, N) >1) \xrightarrow[N \to \infty]{}  \frac{1}{2}.$$
In order to be sure that our observation period remains inside the interval $[0,1]$, the price to pay is to consider a slightly less number of observations, i. e. to take $\alpha \approx 1$ (which means $\alpha<1$ is arbitrary close to $1$). Then, $\tau_{\alpha N} \sim G(\alpha  N , N)$, which can be written as $\tau_{\alpha N} \sim G(M, \tilde{\alpha} M)$ , where $M = \alpha N$ and $\tilde{\alpha} = 1/ \alpha$. We have

\begin{align*}
\mathbb{P} \left( \tau_{\alpha N} > 1 \right) = \int_{1}^{\infty} \dfrac{(\tilde{\alpha} M )^{M}}{\Gamma(M)} x^{M-1} e^{-\tilde{\alpha} M x} dx,
\end{align*}
and by  the change of variable $y = \tilde{\alpha} M x$, the last equation can be written as 
\begin{align*}
\mathbb{P} \left( \tau_{\alpha N} > 1 \right) &= \dfrac{(\tilde{\alpha} M )^{M}}{\Gamma(M)} \int_{\tilde{\alpha} M}^{\infty} \left( \dfrac{y}{\tilde{\alpha} M} \right) ^{M-1} e^{-y} dy \\
&= \dfrac{\tilde{\alpha} M}{\Gamma(M)} \int_{\tilde{\alpha} M}^{\infty} y^{M-1} e^{-y} dy = \dfrac{\tilde{\alpha} M}{\Gamma(M)} \Gamma(M, \tilde{\alpha} M).
\end{align*}

By  the result obtained by \cite{gaut1977} on  the limit behavior of $\Gamma(N, \alpha N)$ for $\tilde{\alpha} > 1$, the last result can be written as
\begin{align*}
\mathbb{P} \left( \tau_{\alpha N} > 1 \right) & \sim \dfrac{\tilde{\alpha} M}{\Gamma(M)} \cdot \dfrac{(\tilde{\alpha} M)^{M} e^{-\tilde{\alpha} M}}{(1 + \tilde{\alpha}) M} = \dfrac{\tilde{\alpha} ^{M+1}  M^{M+1}  e^{-\tilde{\alpha} M}}{(1 + \tilde{\alpha}) M !},
\end{align*}
from Stirling approximation we get, 
\begin{align*}
\mathbb{P} \left( \tau_{\alpha N} > 1 \right) & \sim \dfrac{\tilde{\alpha} ^{M+1}  M^{M+1}  e^{-\tilde{\alpha} M}}{(1 + \tilde{\alpha}) \sqrt{2 \pi M} \left( \frac{M}{e} \right)^{M}} = \dfrac{\tilde{\alpha} ^{M+1}  M^{1 - 1/2}  e^{-\tilde{\alpha} M} e^{M}}{(1 + \tilde{\alpha}) \sqrt{2 \pi} }  \\\
&=\dfrac{\tilde{\alpha} ^{M+1}  M^{1/2}  e^{-M (\tilde{\alpha} +1 )}}{(1 + \tilde{\alpha}) \sqrt{2 \pi} } = \dfrac{\tilde{\alpha} ^{\alpha N+1}  M^{1/2}  e^{-\alpha N (\tilde{\alpha} +1 )}}{(1 + \tilde{\alpha}) \sqrt{2 \pi} } ,
\end{align*}
and recalling that $M = \alpha N$,  the last quantity goes to zero as $N\to \infty$ for $\alpha <1$. 

Then, for every $\alpha \approx 1$ and for every $ \eps >0$ there exists a set $\Omega_{\alpha, \eps}\subset \Omega$ such that $\mathbb{P}(\Omega_{\alpha, \eps} ) >1-\eps $ and for every $\omega \in \Omega _{\alpha, \eps}$ we have $\tau_{[\alpha N]} (\omega) \leq 1$ when $N$ is sufficiently large. We then assume that, in the renewal case,  $\alpha \approx 1$ and we always work on the  space $\Omega_{\alpha, \eps} .$ A similar procedure was considered by \cite{Mi}, Theorem  3.4.1.

Notice that several authors (see e.g. \cite{Vi1}), consider that the model (\ref{1}) is  at times $\tau_{1},..., \tau_{N(1)}$ where $N(1)$ is the last time contained in the interval $[0, 1]$. Instead, we prefer to work on a smaller probability space (but still very close to $\Omega)$  which guaranties that  the entire observation period is contained in the unit interval $[0,1]$. 

\section{Least squares estimator}

Let us fix $0<\alpha \leq 1$ and denote $N_{\alpha }= [\alpha N]$, the number of observations. Consider the model
\begin{equation} \label{modelo}
Y_{\tau_{i+1}} = a \tau_{i+1} + \Delta W_{\tau_{i+1}}, \hskip0.5cm i=0,\ldots N_{\alpha}-1
\end{equation}
where $\Delta W_{\tau_{i+1}}= W_{\tau_{i+1}}-W_{\tau _{i}}$, with $0<\alpha \leq 1$.  Actually, throughout this work we assume  $\alpha =1$ in the jittered sampling case and $\alpha \approx 1$ in the renewal sampling case, see the discussion in Section \ref{sec22}.  In Figures \ref{njs1} and \ref{njs2} we illustrate the behavior of the noise in (\ref{modelo}) at the random times (\ref{js}) and (\ref{rp}) (which appears to be similar to the behavior of the Brownian increment itself). The LSE  for the drift parameter $a$ in the model (\ref{1}) is obtained in a standard way, by minimizing the function $f(a)= \sum_{i=0}^ {N_{\alpha}-1} \left( Y_{\tau_{i+1}} - a \tau_{i+1}\right) ^ {2} $ giving 
\begin{align}
\hat{a}_{N} &=  \dfrac{\displaystyle  \sum_{i=0}^{N_{\alpha}-1} \tau_{i+1} Y_{\tau_{i+1}}}{\displaystyle  \sum_{i=0}^{N_{\alpha}-1} \tau_{i+1}^2}  \label{a-jt} 
\end{align}
for both jittered sampling (JS) and renewal sampling  (RS) cases. From (\ref{modelo}), (\ref{a-jt}) we immediately have, 

\begin{align}
\hat{a}_{N} - a  &=  \dfrac{\displaystyle \dfrac{1}{N} \sum_{i=0}^{N_{\alpha}-1} \tau_{i+1} \Delta W_{\tau_{i+1}}}{\displaystyle \dfrac{1}{N} \sum_{i=0}^{N_{\alpha}-1} \tau_{i+1}^2} := \dfrac{A_N}{D_N}, \quad
\label{DN-js}
\end{align} 
for every $N\geq 1$,
\begin{equation}
\label{AD}
A_{N}=  \sum_{i=0}^{N_{\alpha}-1} \tau_{i+1} Y_{\tau_{i+1}} \mbox{ and } D_{N}=\dfrac{1}{N} \sum_{i=0}^{N_{\alpha}-1} \tau_{i+1}^2.
\end{equation}

Our purpose is to analyze the asymptotic properties of the LSE (\ref{a-jt}), in particular its asymptotic normality in distribution. The denominator of the expression (\ref{DN-js}) has been already studied by \cite{araya2019}. Let us recall their results (see Lemma 3.2 of \cite{araya2019}). 

\begin{prop}\label{pp1}
Let $D_{N}$ be given by (\ref{AD}). Then $D_{N}$ converges almost surely, as $N\to \infty$ to $\frac{\alpha^{3}}{3}$. 

\end{prop}

Actually, the result of \cite{araya2019}  has been obtained for $\alpha=1$, but after inspecting the proof, it is clear that the same arguments holds for every $\alpha \in (0, 1)$. Therefore, in order to obtain the asymptotic behavior of the LSE, we need to analyse the sequence $A_{N}$ in (\ref{DN-js}).  A first step in this direction is to evaluate the $L^ {2}(\Omega)$- norm of $A_{N}$ when $N$ is large.

\begin{lemma}\label{ll1}
Let $A_{N}$ given by \eqref{DN-js}, for every $N \geq 1$, then either if $\tau$ is defined as (\ref{js}) or (\ref{rp}), it holds
\begin{equation*}
\mathbb{E} \left| N A_{N} \right| ^{2} \xrightarrow[N \to \infty]{} \frac{1}{3}\alpha ^ {3}.
\end{equation*}
\end{lemma}

Although the  above result is the same in the JS and RS cases, the proof is different. While in the JS case, the limit if given by the  ``deterministic part'' of the times (\ref{js}), in the RS case there is no deterministic part and both summands in (\ref{E2-RS}) contribute to the limit. The prroofs can be found in the appendix.

\section{Limit distribution of the LSE}
Lemma \ref{ll1} shows that the sequence $(A_{N}) _{N\geq 1}$ given by (\ref{AD}) converges in $L^ {2}(\Omega)$ to zero as $N\to \infty$. We can also  show that  $A_{N}$ converges to zero in $L^ {p}(\Omega)$ for every $p\geq 2$ and by a Borel-Cantelli argument, we get its almost sure convergence to zero. Indeed, via conditioning on $\nu$, 
$$\mathbb{E} \vert A_{N}\vert ^ {p} = \mathbb{E} \left[ \mathbb{E} \left[ \vert A_{N}\vert ^ {p}| \nu \right] \right]=\mathbb{E} \left[ g( \tau _{1},\ldots, \tau_{N})\right]$$
with, for $x_{1}< x_{2}<\ldots x_{N}$, 

$$g(x_{1},..., x_{N})= \mathbb{E} \left| \sum_{i=0} ^ {N-1} x_{i+1} (W_{x_{i+1}}-W_{x_{i}} ) \right| ^ {p} \leq C_{p} \left(\mathbb{E} \left| \sum_{i=0} ^ {N-1} x_{i+1} (W_{x_{i+1}}-W_{x_{i}} ) \right| ^ {2} \right) ^ {\frac{p}{2}}.$$ This implies, together with Lemma \ref{ll1} 
$$\mathbb{E} \vert A_{N}\vert ^ {p} \leq C_{p} N ^ {-2p}$$
and thus for every $\gamma >0$  and for every $p$ integer such  $2p-\gamma>1$, 
$$\sum_{N\geq 1} \mathbb{P}( A_{N} > N ^ {-\gamma} ) \leq \sum_{N\geq 1} N ^ {\gamma-2p} <\infty.$$
Then, via Proposition \ref{pp1}, we obtain the consistency of the LSE  (\ref{a-jt}).

Let us now study the asymptotic limit in distribution of $A_{N}$. To this end, we need to study the sequence $(Q_{N})_{N\geq 1}$ defined, for every $N\geq 1$, by
\begin{align} \label{qn}
Q_{N} &=  \sum_{j=0}^{N_{\alpha}-1} \tau_{j+1}^{2} \left( \tau_{j+1} - \tau_{j} \right). 
\end{align} 
This plays the role of the ``bracket'' of $A_{N}$. Before, let us introduce some notation: If $\tau_{j}$ is given by (\ref{js}), then we can write them as 
\begin{equation}\label{3d-1}
\tau_{j} = \frac{j}{N} + \frac{X_{j}}{N}, 
\end{equation}
where $X_{j}, j=1,..,.$ are independent random variables and $X_{j}$ follows a symmetric probability distribution with support in  $\left[ -\frac{1}{2} , \frac{1}{2} \right]$ and with density denoted by $g$. If $\tau_{j}$ is given by (\ref{rp}), then 
\begin{equation}\label{3d-2}
\tau_{j} = \frac{1}{N} \sum_{i=1}^{j} X_{i}, 
\end{equation}
where $X_{i}$ follows the Gamma law $G(1, 1)$. Moreover, for every $i\geq 1$, $X_{i+1}- X_{i}$ is independent of $X_{l}$ if $l\leq i$. 

\begin{prop} \label{L2-conv}
Let $Q_{N}$ be given by (\ref{qn}). Then
 $$\lim_{N \to \infty} \mathbb{E} \left[ \left( Q_{N} - \frac{1}{3} \right)^{2} \right]   \xrightarrow[N \to \infty]{} 0.$$
\end{prop}

We now give the limit in distribution of the sequence $A_{N}$.

\begin{prop} \label{normalidad}
Let $A_{N}$ given in \eqref{DN-js}, then the following convergence in distribution holds
\begin{equation*}
N A_{N} \xrightarrow[N \to \infty ^{+}]{\mathcal{L}} N \left( 0, \frac{\alpha ^ {3}}{3} \right),
\end{equation*}
either if $\tau _{i}$ are defined  either by  \eqref{js} or by  \eqref{rp}.
\end{prop}
{\bf Proof: }   We  analyse the asymptotic distribution of the characteristic function of $NA_{N}$, denoted $\varphi_{NA_{N}}$ in the sequel.  Via conditioning, with $X_{j}$ given by (\ref{3d-1})  or by (\ref{3d-2}), with $N'_{\alpha}= N-1$ in the JS case and $N'_{\alpha} = N_{\alpha} $ in the RP case, 
\begin{align}
\varphi_{N A_{N}}(t) = \mathbb{E} \left[ e^{it N A_{N}} \right]
&= \mathbb{E} \left[ e^{it \sum_{j=0}^{N'_{\alpha}-1} \left( \frac{j+1}{N} + \frac{X_{j+1}}{N} \right) \left( W_{\frac{j+1}{N} + \frac{X_{j+1}}{N}} - W_{\frac{j}{N} + \frac{X_{j}}{N}} \right) } \right] \nonumber \\
&= \mathbb{E} \left. \left[ \mathbb{E} \left[ e^{it \sum_{j=0}^{N'_{\alpha}-1} \left( \frac{j+1}{N} + \frac{X_{j+1}}{N} \right) \left( W_{\frac{j+1}{N} + \frac{X_{j+1}}{N}} - W_{\frac{j}{N} + \frac{X_{j}}{N} }\right) } \right| X \right] \right] \nonumber \\
&= \mathbb{E} \left[ e^{-\frac{t^{2}}{2} \sum_{j=0}^{N'_{\alpha}-1} \left( \frac{j+1}{N} + \frac{X_{j+1}}{N} \right)^{2} \left( \frac{X_{j+1}}{N} - \frac{X_{j}}{N} + \frac{1}{N}  \right) } \right]=\mathbb{E} \left[ e^{-\frac{t^{2}}{2} Q_{N} } \right]   \label{lim-js-bm}
\end{align}
with $Q_{N}$ from (\ref{qn}). Now, by Propostion \ref{L2-conv}, the sequence $Q_{N}$ converges, in $L^ {2}(\Omega)$, thus in probability, to $\frac{\alpha ^ {3}}{3}$. By Dominated convergence theorem, for every $t\in \mathbb{R}$, 

$$ \lim_{N \to \infty}\varphi_{N A_{N}}(t)=  \mathbb{E} \left[  e^{-\frac{t^{2}}{2} \frac{\alpha^{3}}{3}} \right]$$
and this gives the conclusion. \qed

By Propositions \ref{pp1} and \ref{normalidad}, we immediately obtain the asymptotic normality of the LSE. We denote by $\xrightarrow{\mathcal{L}}$ the convergence in law. 

\begin{theorem} \label{an_dn}
Consider the LSE $\hat{a}_{N}$ given by (\ref{a-jt}). Then 
$$N( \hat{a}_{N} - a) \xrightarrow[N \to \infty ^{+}]{\mathcal{L}}N(0,\frac{3}{\alpha ^{3}}). $$
\end{theorem}

\begin{remark}
Let us give some heuristics that explain the convergence in law of the sequence $(NA_{N}) _{N\geq 1}$.  Consider the JS case and $\alpha =1$.  Then we can write 
\begin{eqnarray*}
NA_{N}=\sum_{i=0} ^{N-2} \tau_{i+1} (W_{\tau_{i+1}}- W_{\tau_{i}})= \int_{0}^{1} H_{N} (s) dW_{s}
\end{eqnarray*}
with $H_{N}(s)= \sum_{i=0} ^{ N-2} \tau_{i+1} 1_{(\tau_{i}, \tau_{i+1}]}(s)$, for $s\in [0,1]$.  Intuitively, $H_{N}(t)$
 converges to $t$ in $L ^{2}([0,1]\times \Omega)$ since $\vert \tau_{i}-\frac{i}{N} \vert \leq \frac{1}{N}$. Therefore 
$NA_{N}$ would converge in $L^{2}(\Omega)$, as $N\to \infty$,  to $\int_{0}^{1} sdW_{s}$ whose law in $N(0, \frac{1}{3})$. 

\end{remark}

Let us finish this theoretical part with some comment on the distance between the law of the sequence $(NA_{N})_{N\geq 1}$  and its limit. Recall that the distance between the laws of two random variables $X$ and $Y$ is defined as
$$d (X, Y) =\sup_{h\in\mathcal{A}} \left| \mathbb{E} h(X)- \mathbb{E}h(Y) \right| $$
where $ \mathcal{A}$ is a class of functions (its choice  defines specific distances, such as Kolmogorov, total variation or Wasserstein or other distances). Let $h$ be a function such that the all expectations below exist. Consider for simplicity $\alpha=1$ as in the JS case. Then, by taking the conditional expectation as in the proof of Lemma \ref{ll1}, we obtain 
$$\mathbb{E} h(NA_{N}) = \mathbb{E} h \left( \sqrt{\mathbb{E}(NA_{N}) ^ {2}} Z\right)$$
with $Z\sim N(0,1)$. This implies that, for $N$ large (see Proposition 3.6.1 of \cite{NP})
\begin{equation}
\label{6d-1}d  \left(NA_{N}, \sqrt{\frac{1}{3}}Z\right) \leq C \left| \mathbb{\E}(NA_{N})^ {2} -\frac{1}{3} \right| \leq C \frac{1}{N}
\end{equation}
where the last bound can be obtained easily from the proof of Lemma \ref{ll1}. A similar bound as in (\ref{6d-1}) can be obtained when we replace $NA_{N}$ by $\hat{a}_{N}-a$ with $\hat{a}_{N} $ given by (\ref{a-jt}).

\section{Simulation study}

In this section we consider the different problems that appear when studying the limit distribution of the least square estimator $\hat{a}_{N}$, properly normalized. First, we show the behavior of the increment of the standard Brownian motion, considering both type of random times, for different values of $N$. Next, we illustrate how the number of times, in which $ \tau_{\alpha N} > 1$, changes for different values of $\alpha$. In addition, we show by simulation the  convergence of $Q_N$ as we proved in  Proposition \ref{L2-conv} and \ref{normalidad}, which are necessary to show Theorem \ref{an_dn}. Finally, we define the error of our estimation and the corresponding simulation result. \\

We have simulated the observations $Y_{\tau_{1}}, \dots, Y_{\tau_{N_{\alpha}}}$, considering $N=5000$ and $\alpha = 1$ for jittered sampling  and $\alpha = 0.98$ for renewal sampling, and we have repeated this procedure 10000 times in order to obtain the corresponding tables and histograms

\subsection*{Increment of standard Brownian motion under observations sampled at random times:}
Let us recall that the standard Brownian motion is an adapted process defined in some probability space $(\Omega, \mathcal{F}, \mathbb{P})$, $W_{0} = 0$, it has independent and stationary increments which follow a normal distribution, i.e. $W_{t} - W_{s} \sim N(0, t-s)$. In the following figures, it is possible to notice the behavior of the increment of the standard Brownian motion for different values of N and the different types of random times defined above.
\begin{figure}[h!]
\centering
\begin{subfigure}{.40\textwidth}
  \centering
  \includegraphics[width=.95\linewidth]{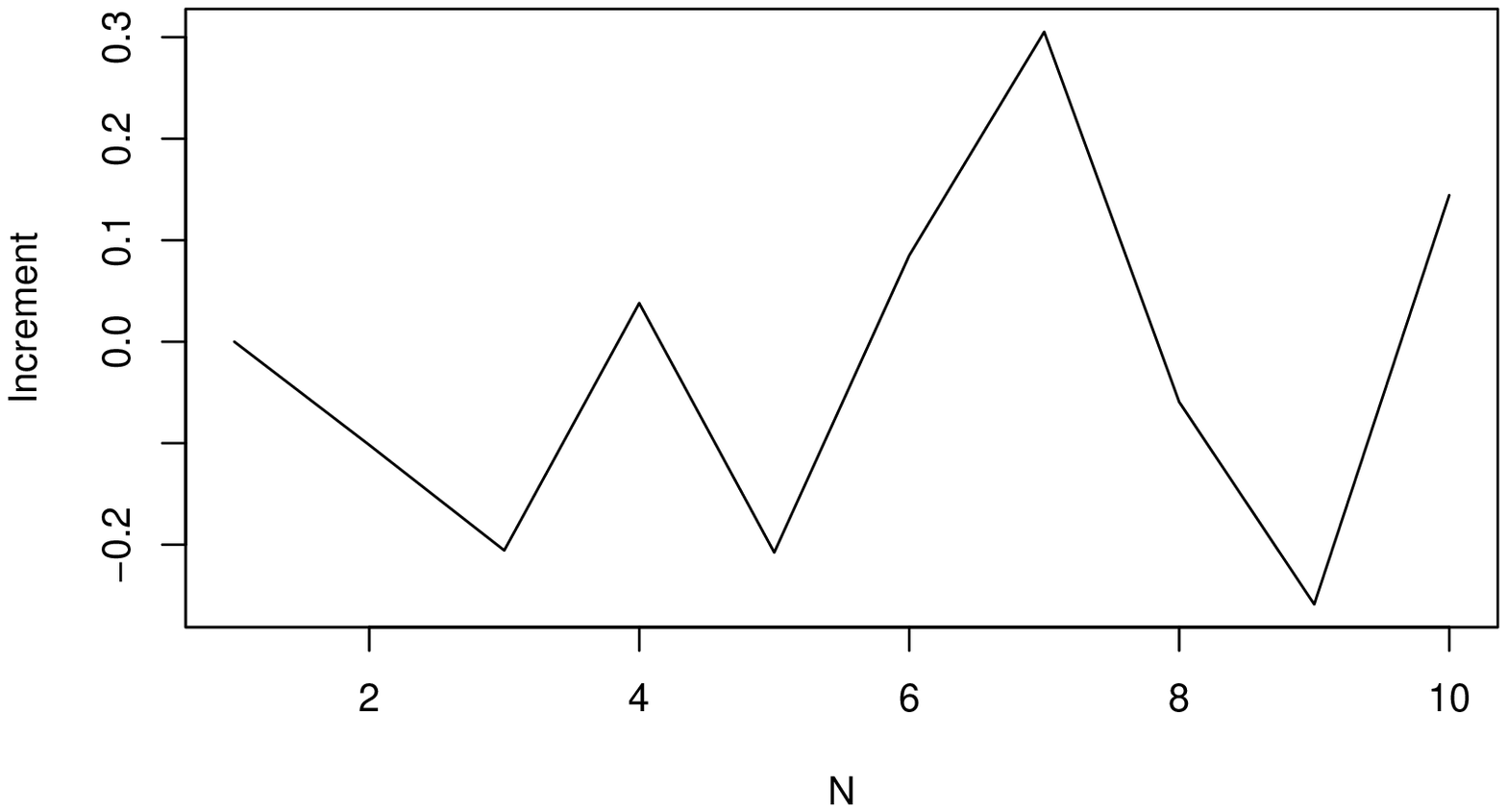}
  \label{n10js}
\end{subfigure}%
\begin{subfigure}{.40\textwidth}
  \centering
  \includegraphics[width=.95\linewidth]{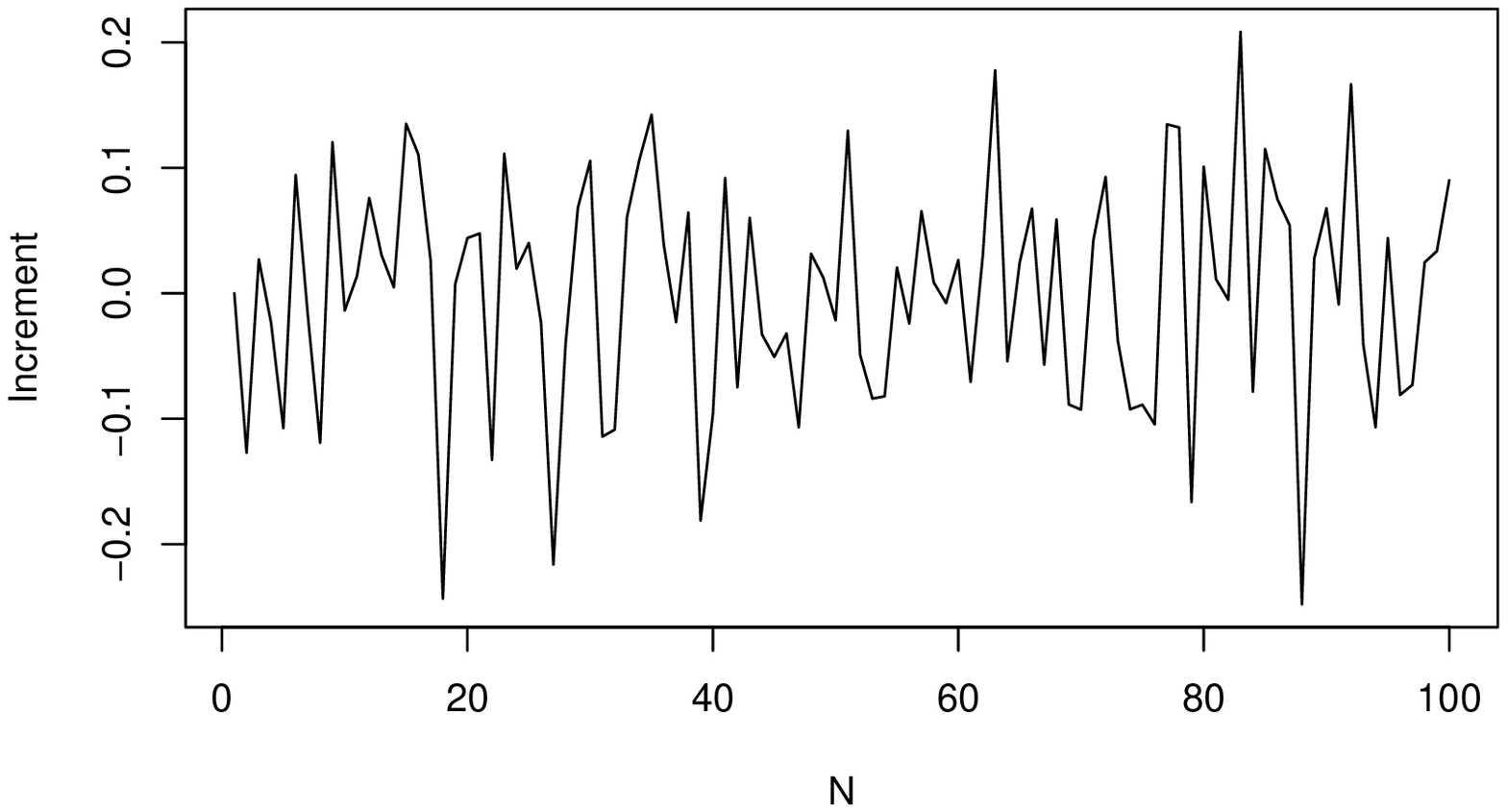}
  \label{n100js}
\end{subfigure}
\begin{subfigure}{.42\textwidth}
  \centering
  \includegraphics[width=.95\linewidth]{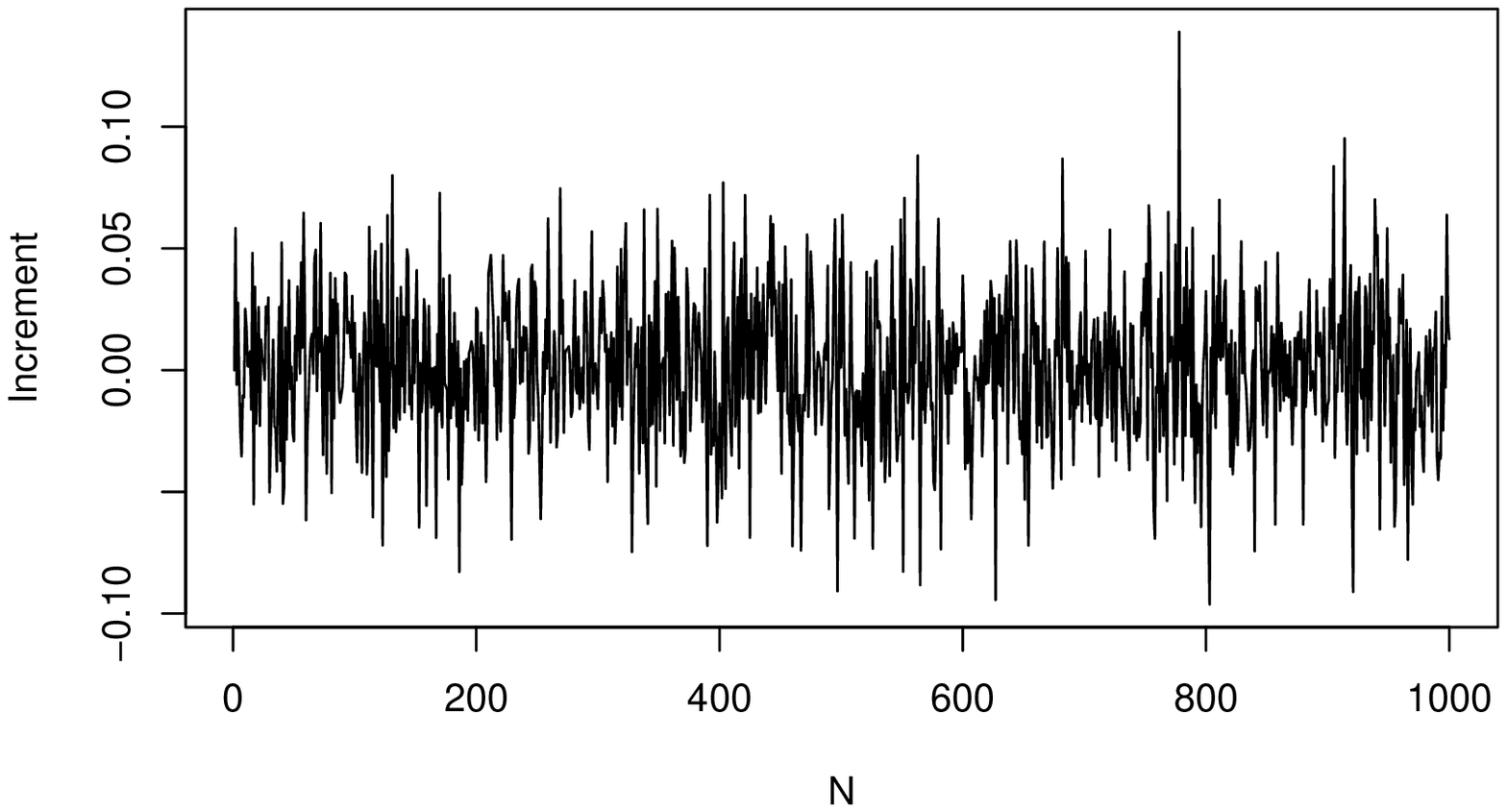}
  \label{n1000js}
\end{subfigure}
\caption{Behavior of the increment of standard Brownian motion under Jittered sampling observations with N=10, 100 and 1000 respectively.}
\label{njs1}
\end{figure}

\begin{figure}[h!]
\centering
\begin{subfigure}{.40\textwidth}
  \centering
  \includegraphics[width=.95\linewidth]{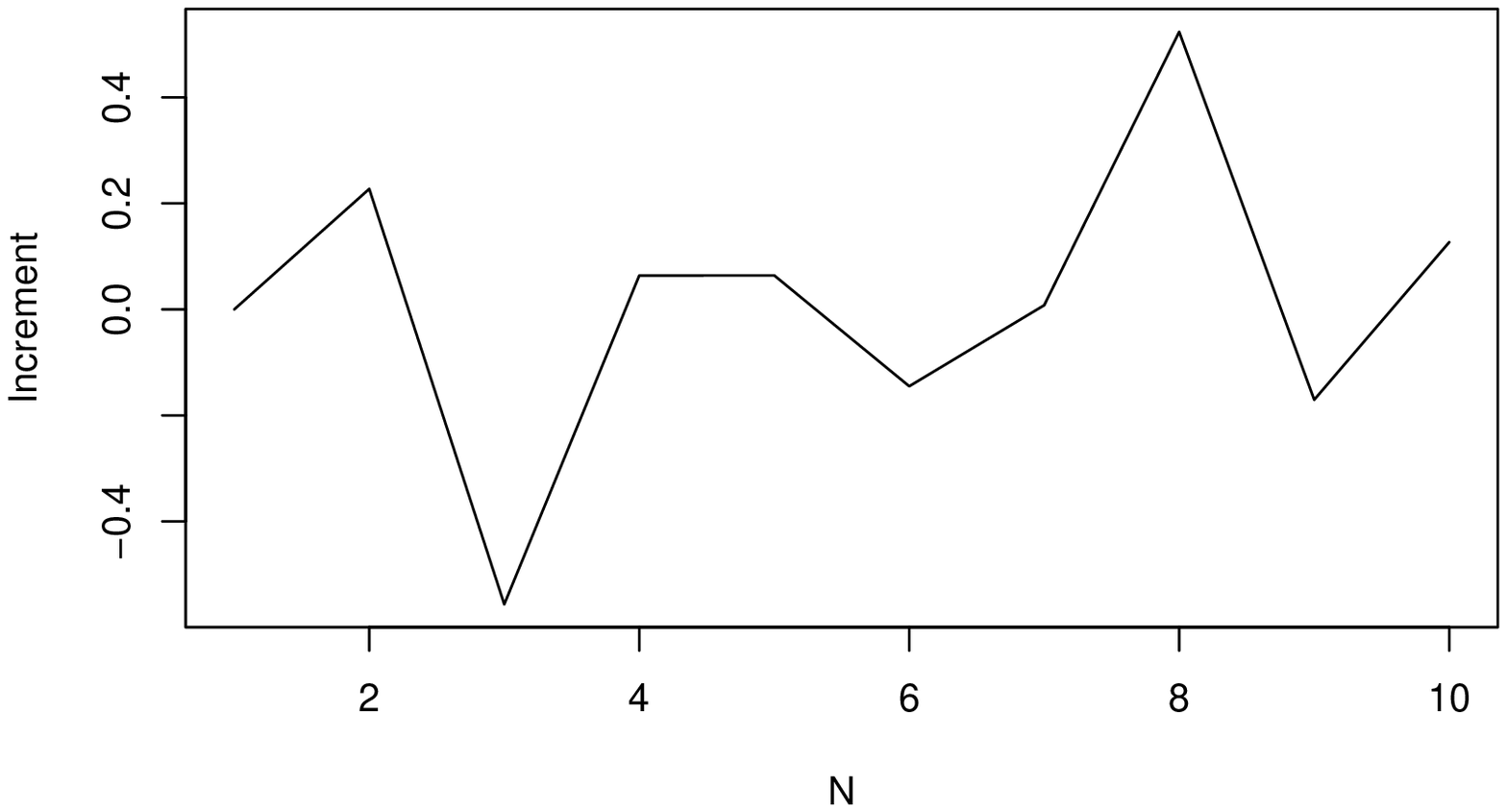}
  \label{n10rp}
\end{subfigure}%
\begin{subfigure}{.40\textwidth}
  \centering
  \includegraphics[width=.95\linewidth]{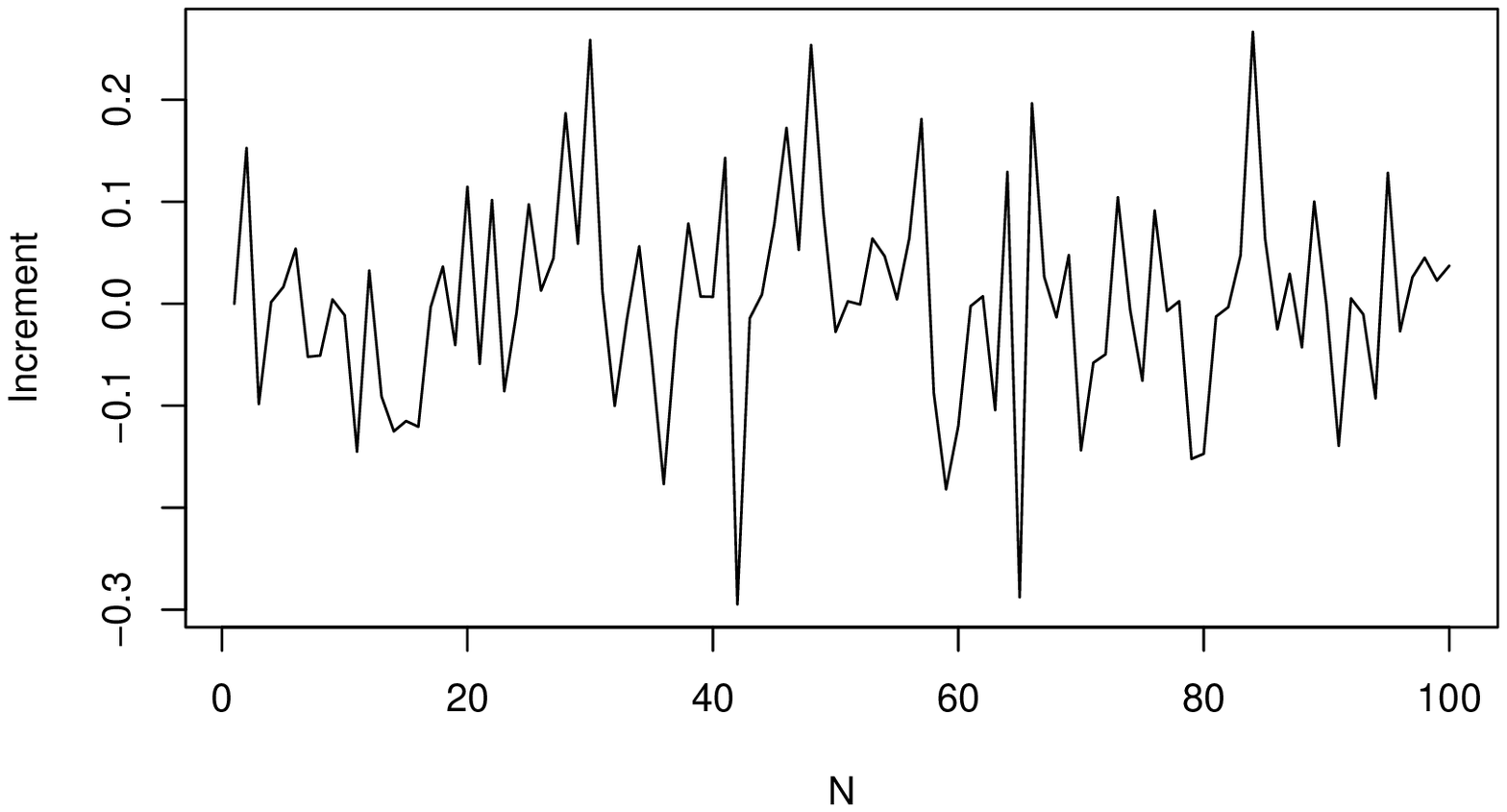}
  \label{n100rp}
\end{subfigure}
\begin{subfigure}{.42\textwidth}
  \centering
  \includegraphics[width=.95\linewidth]{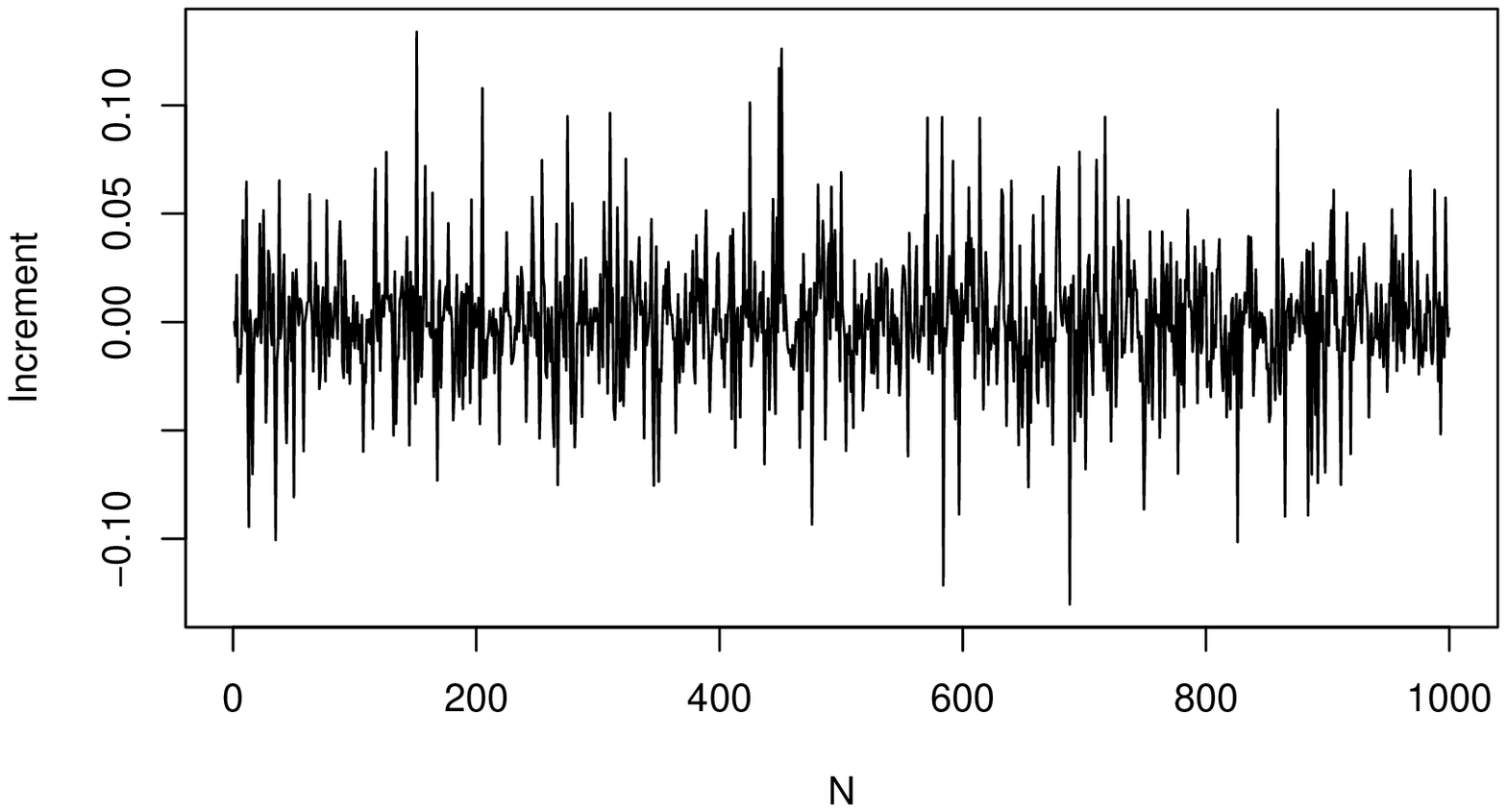}
  \label{n1000rp}
\end{subfigure}
\caption{Behavior of the increment of standard Brownian motion under Renewal sampling observations with N=10, 100 and 1000 respectively.}
\label{njs2}
\end{figure}

\subsection*{Number of times  $\boldsymbol{\tau_{\alpha N} >1}$:}
Under renewal sampling observations, for different values of $\alpha$ and $N$, we have the following experimental results concerning the number of times when the last observations is bigger that 1 (see the discussion in Section 2.2) 
\begin{table}[h!] 
\centering
\begin{tabular}{cccc}
\hline \hline
                             & $N=100$ & $N=1000$ & $N=10000$ \\ \hline
$\alpha = 0.99$ & 45 & 361 & 1624 \\
$\alpha = 0.98$ & 34 & 241 & 225 \\
$\alpha = 0.97$ & 29 & 164 & 17 \\ 
$\alpha = 0.96$ & 27 & 102 & 0 \\
$\alpha = 0.95$ & 24 & 50 & 0 \\
\hline \hline
\end{tabular}
\caption{Number of times $\tau_{\alpha N} >1$ for different values of $\alpha$ and $N$} \label{alpha_tau}
\end{table}

\subsection*{Convergence of  $\boldsymbol{ Q_N}$, Proposition \ref{L2-conv}}
As we shown in Proposition \ref{L2-conv}, $Q_{N}$ converges in mean square and also in probability to $1/3$ as $N$ goes to infinity. After simulating the value of $Q_{N}$ for different values of $N = 2, \dots , 10000$,  we obtain the shape of the convergence. 
In the  following Figure \ref{conv-plot}, the straight line in red, represents the value $Q_{N} = 1/3$. 

\begin{figure}[h!]
\centering
\begin{subfigure}{.5\textwidth}
  \centering
  \includegraphics[width=.9\linewidth]{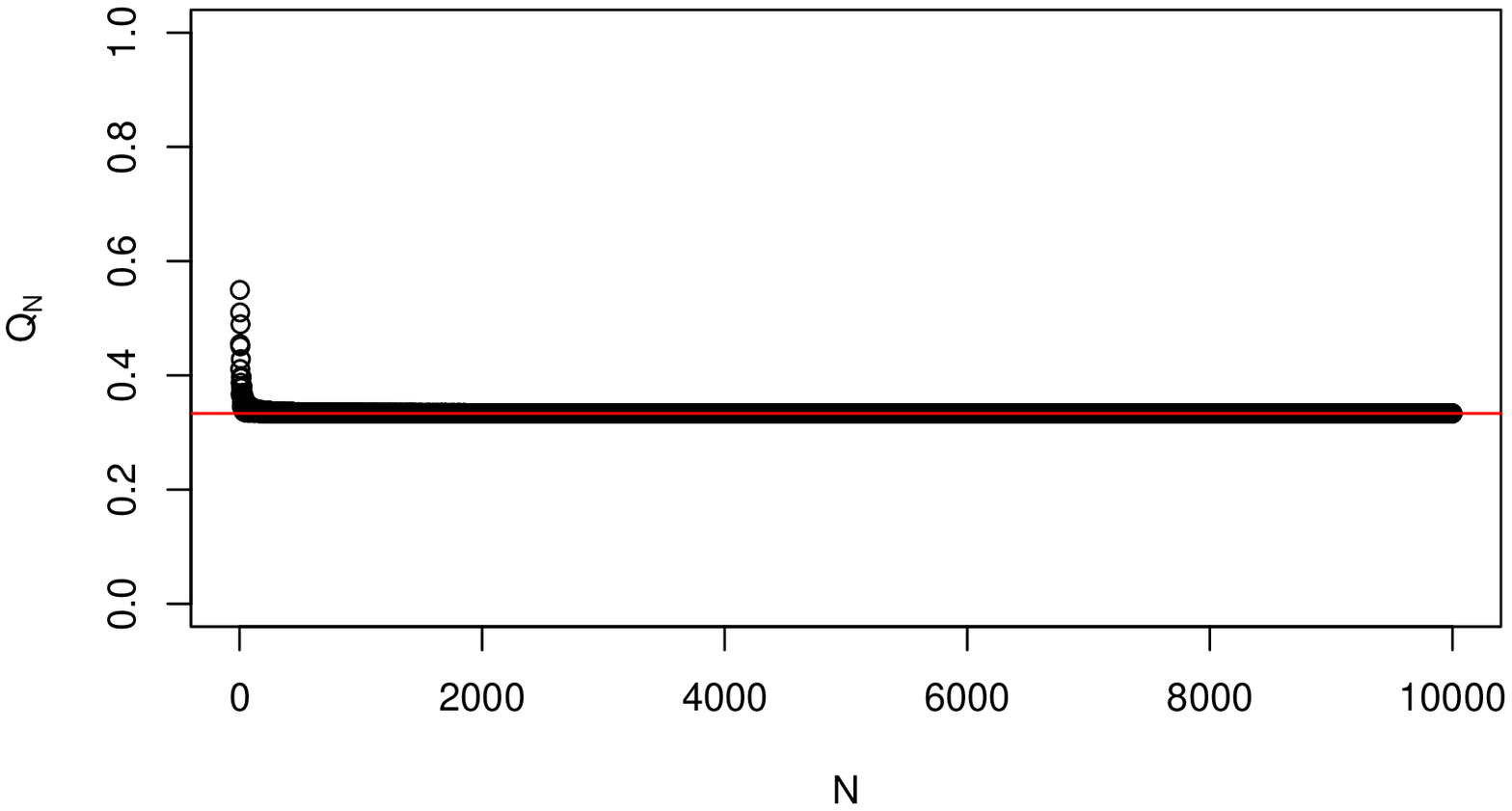}
  \label{ejs}
\end{subfigure}%
\begin{subfigure}{.5\textwidth}
  \centering
  \includegraphics[width=.9\linewidth]{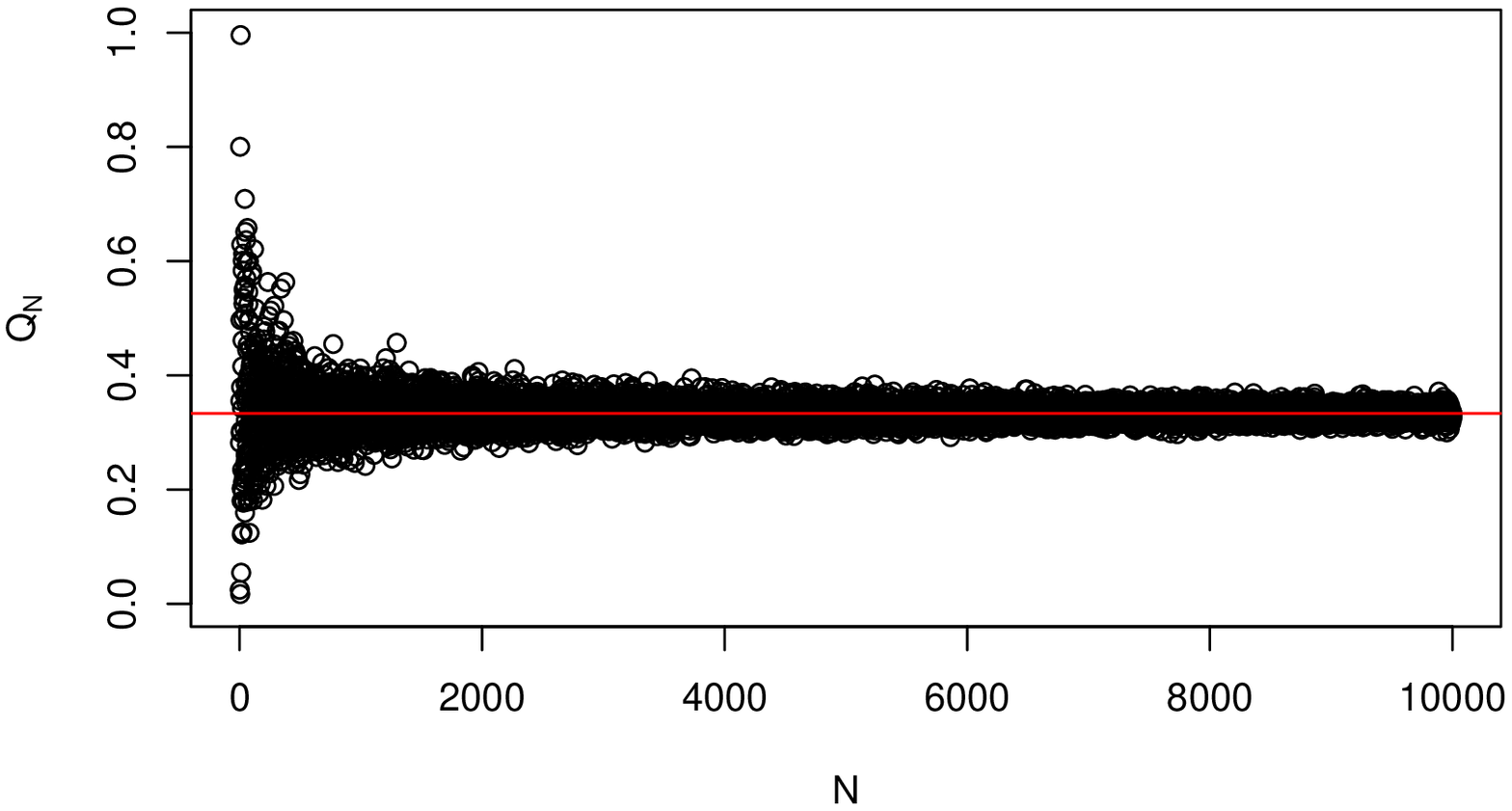}
  \label{erp}
\end{subfigure}
\caption{Convergence of the sequence $Q_N$ under different random times. Left: Jittered Sampling, Right: Renewal Sampling.}
\label{conv-plot}
\end{figure}

\pagebreak
\subsection*{Convergence of  $\boldsymbol{ Q_N}$, Proposition \ref{normalidad}}
 The results obtained, for both type of random times, are the following
\begin{figure}[h!]
\centering
\begin{subfigure}{.5\textwidth}
  \centering
  \includegraphics[width=.9\linewidth]{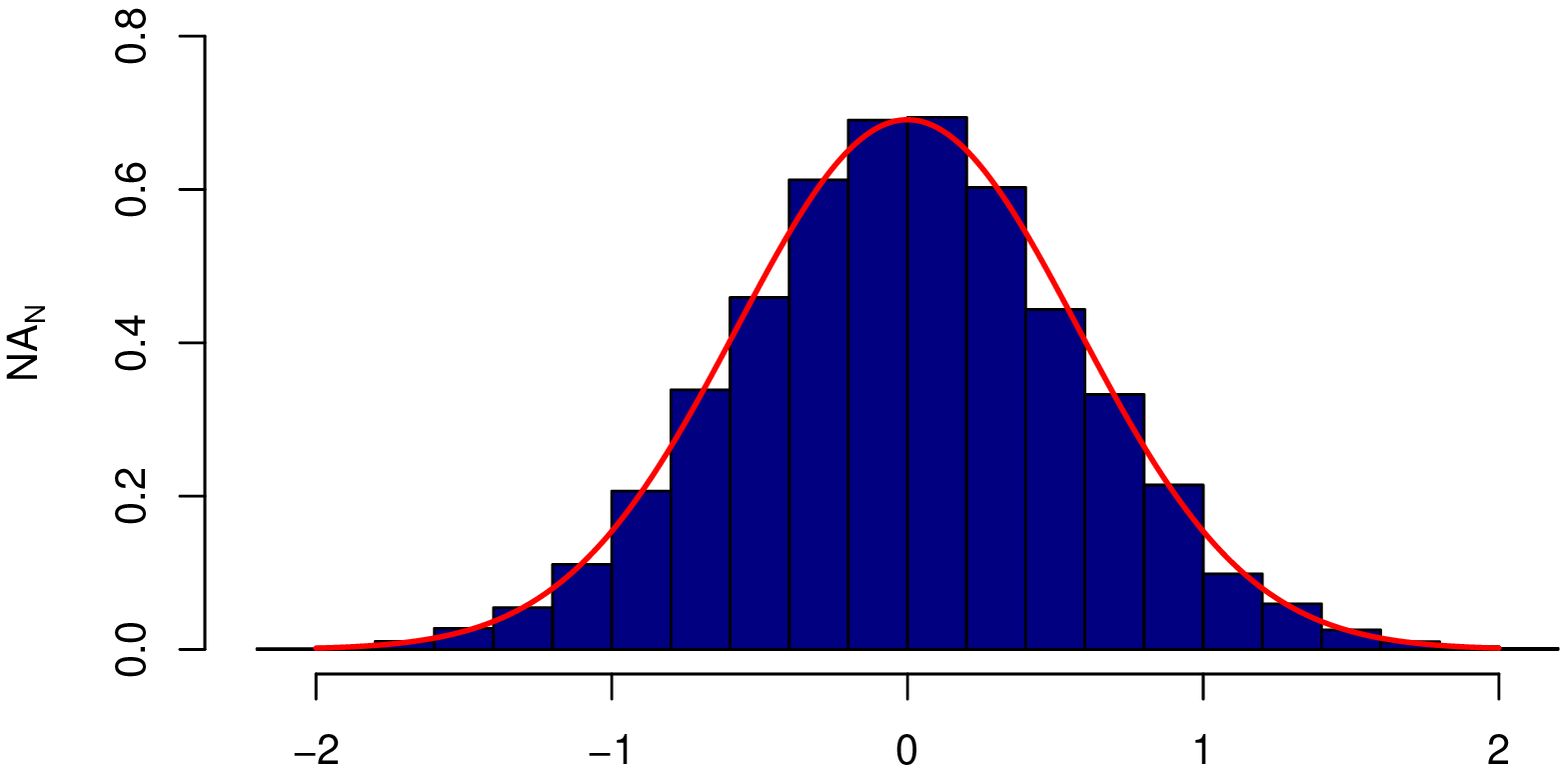}
  \label{ejs1}
\end{subfigure}%
\begin{subfigure}{.5\textwidth}
  \centering
  \includegraphics[width=.9\linewidth]{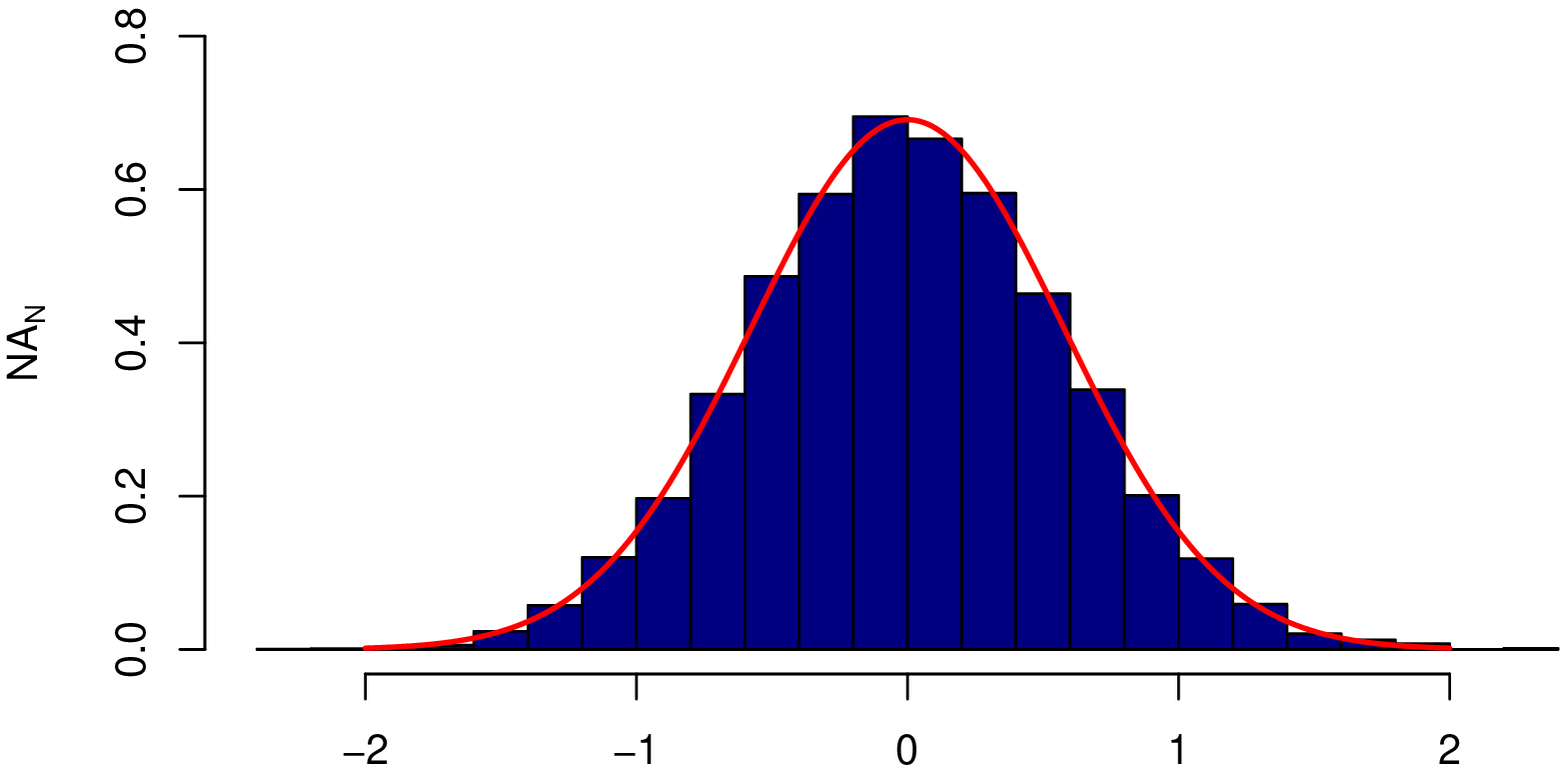}
  \label{erp1}
\end{subfigure}
\caption{Histograms of the sequence $N A_{N}$ under different random times. Left: Jittered Sampling, Right: Renewal Sampling}
\label{hist-plot}
\end{figure}

\begin{table}[h!] 
\centering
\begin{tabular}{ccc}
\hline \hline 
$N A_{N}$& Mean  & Variance \\ \hline
Jittered Sampling & -0.003562114  & 0.3328107 \\
Renewal Sampling  & 0.002266849 & 0.3348508 \\
\hline \hline 
\end{tabular}
\caption{Mean and Variance of the sequence $N A_{N}$} \label{hist-tab}
\end{table}

\newpage
\subsection*{Convergence of $\boldsymbol{A_N}$: Theorem \ref{an_dn}}
\begin{figure}[h!]
\centering
\begin{subfigure}{.5\textwidth}
  \centering
  \includegraphics[width=.9\linewidth]{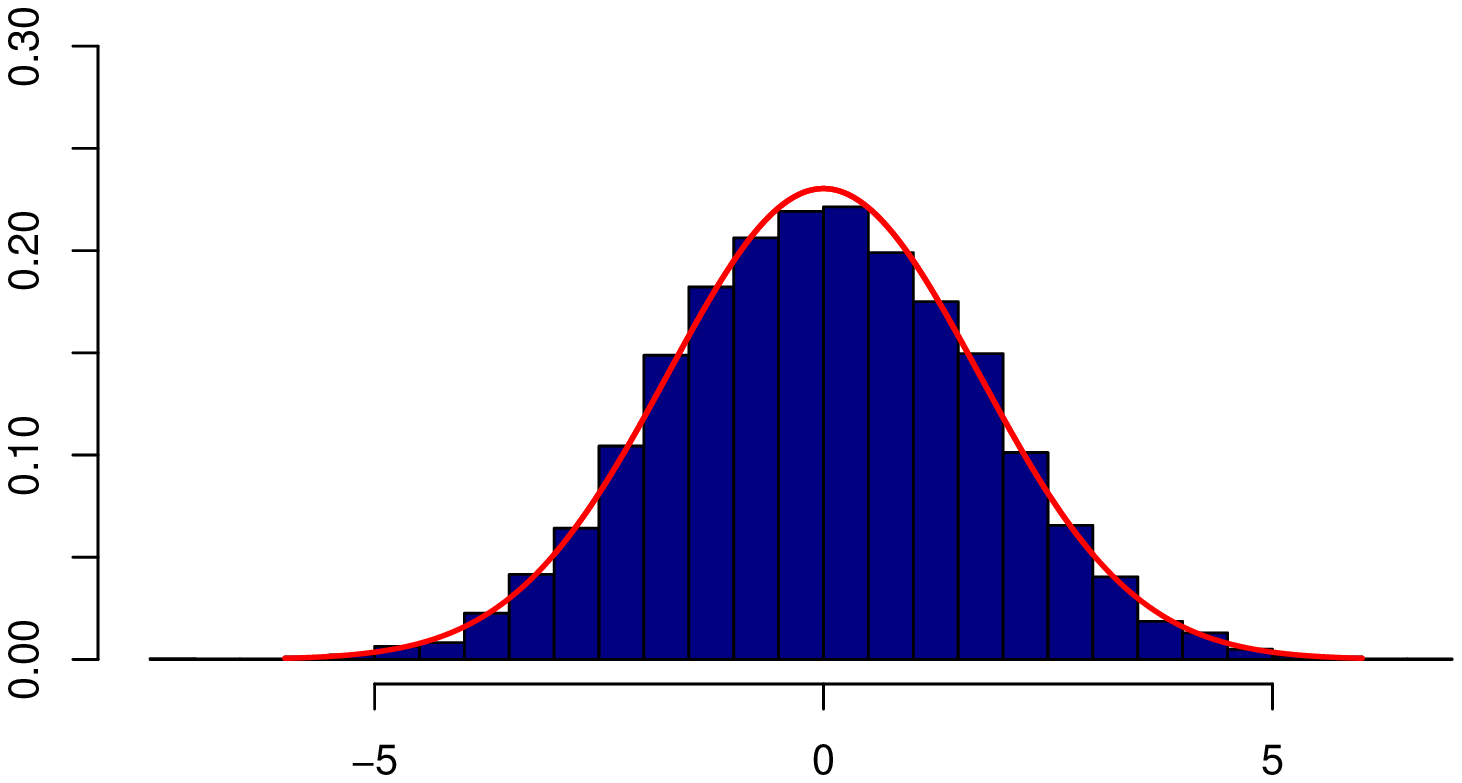}
  \label{teo_js_hist}
\end{subfigure}%
\begin{subfigure}{.5\textwidth}
  \centering
  \includegraphics[width=.9\linewidth]{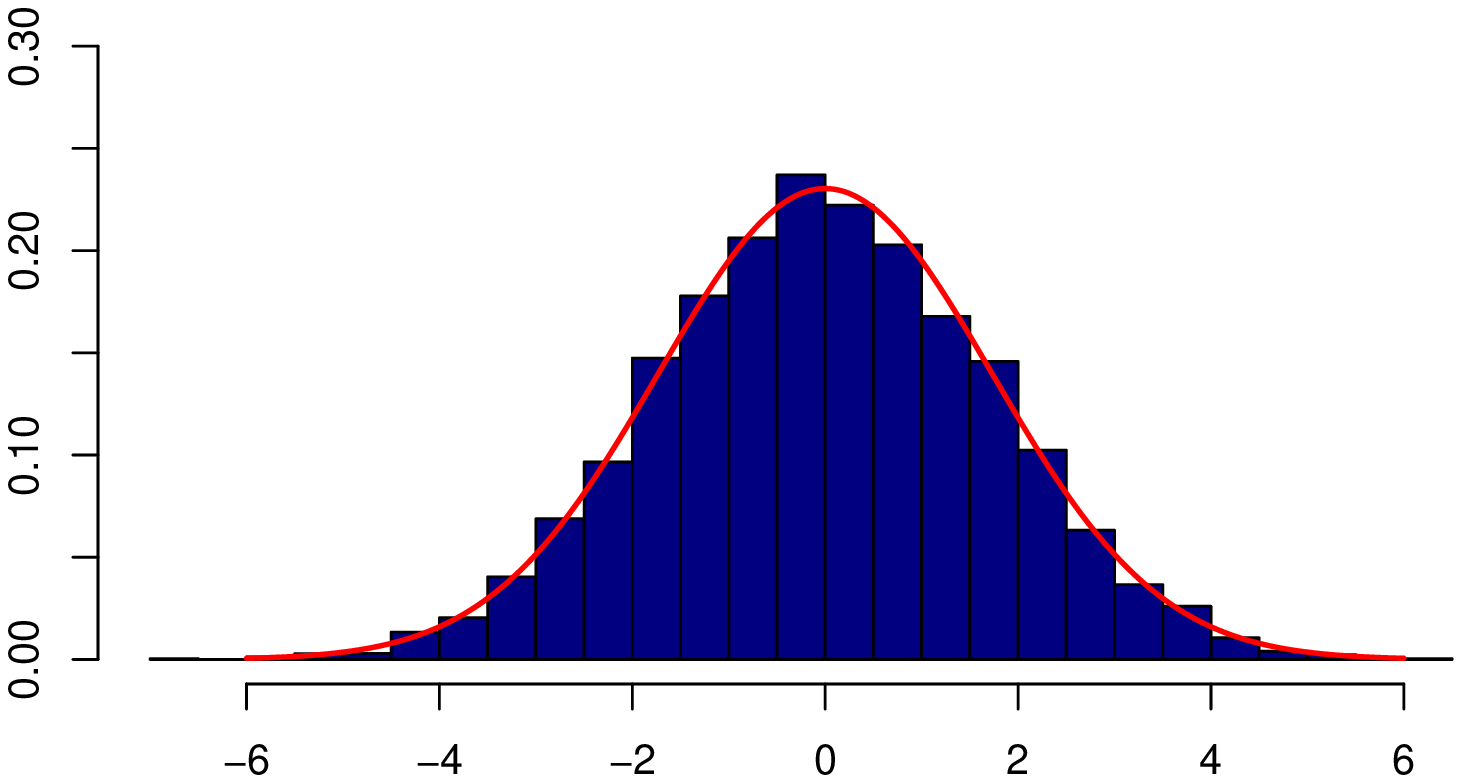}
  \label{teo_rp_hist}
\end{subfigure}
\caption{Histograms of the sequence $N \hat{a}_{N}$  under different random times. Left: Jittered Sampling, Right: Renewal Sampling}
\label{teohist}
\end{figure}

\begin{table}[h!]
\centering
\begin{tabular}{ccc}
\hline \hline 
Random Time & Mean & Variance \\ \hline
Jittered Sampling & -0.001518635 & 3.016777 \\
Renewal Sampling & -0.01345249 & 2.982566 \\ 
\hline \hline 
\end{tabular}
\caption{Mean and variance of the sequence $N \hat{a}_{N}$} \label{tab-teo}
\end{table}

\subsection*{Estimation Error}
As we shown previously, the sequence $N A_{N}$ converges, in law for $N$ large enough, to a normal distribution with mean $\mu =0$ and variance $\sigma^{2} = 1/3$, so we define the estimation error, for both type of random times, as $\varepsilon = N A_{N} - N(0,1/3)$. First, we present the behavior of the error for different values of $N$. The corresponding results, for both type of random times, are the followings

\begin{figure}[h!]
\centering
\begin{subfigure}{.5\textwidth}
  \centering
  \includegraphics[width=.9\linewidth]{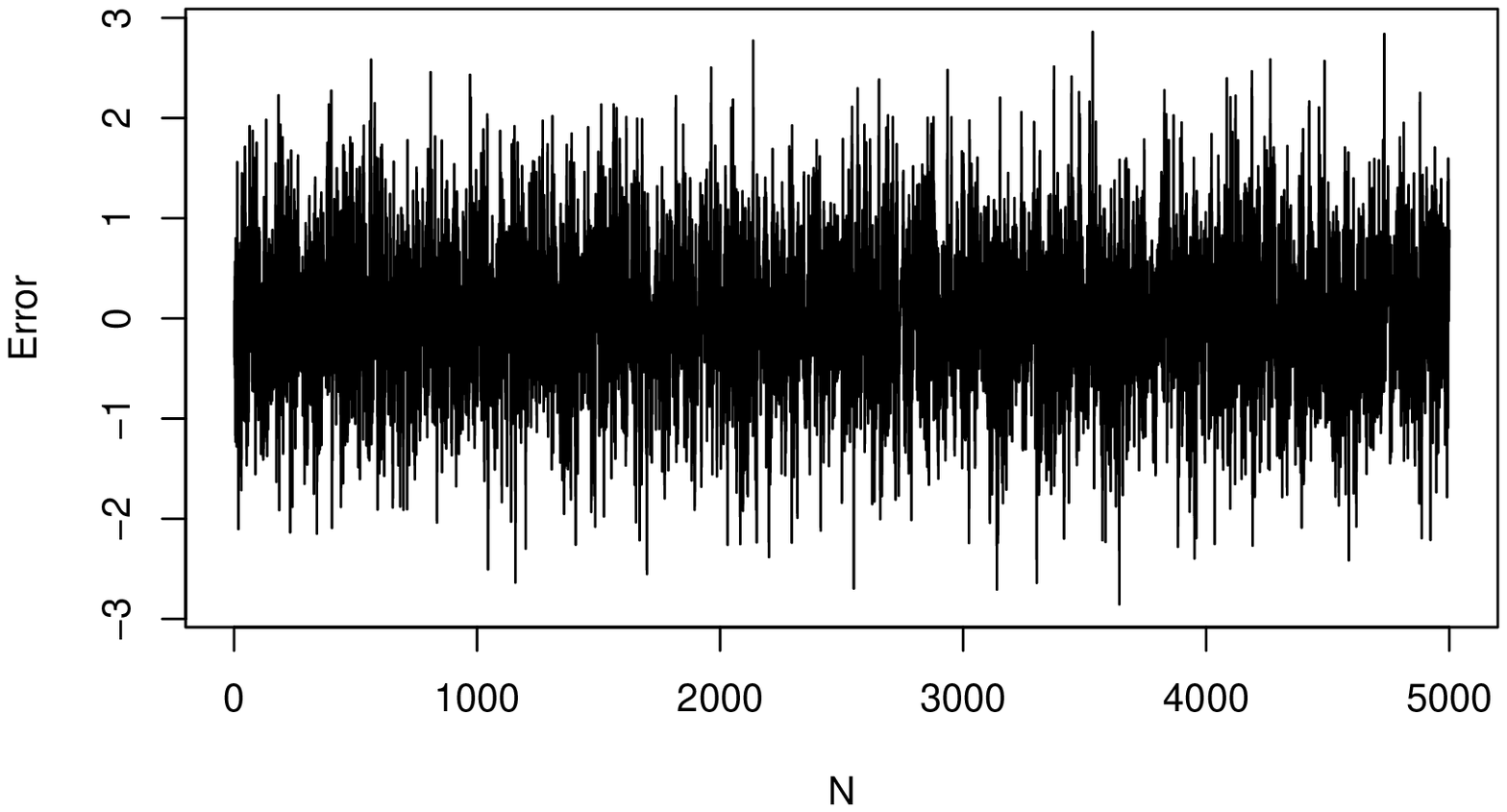}
  \label{ejs2}
\end{subfigure}%
\begin{subfigure}{.5\textwidth}
  \centering
  \includegraphics[width=.9\linewidth]{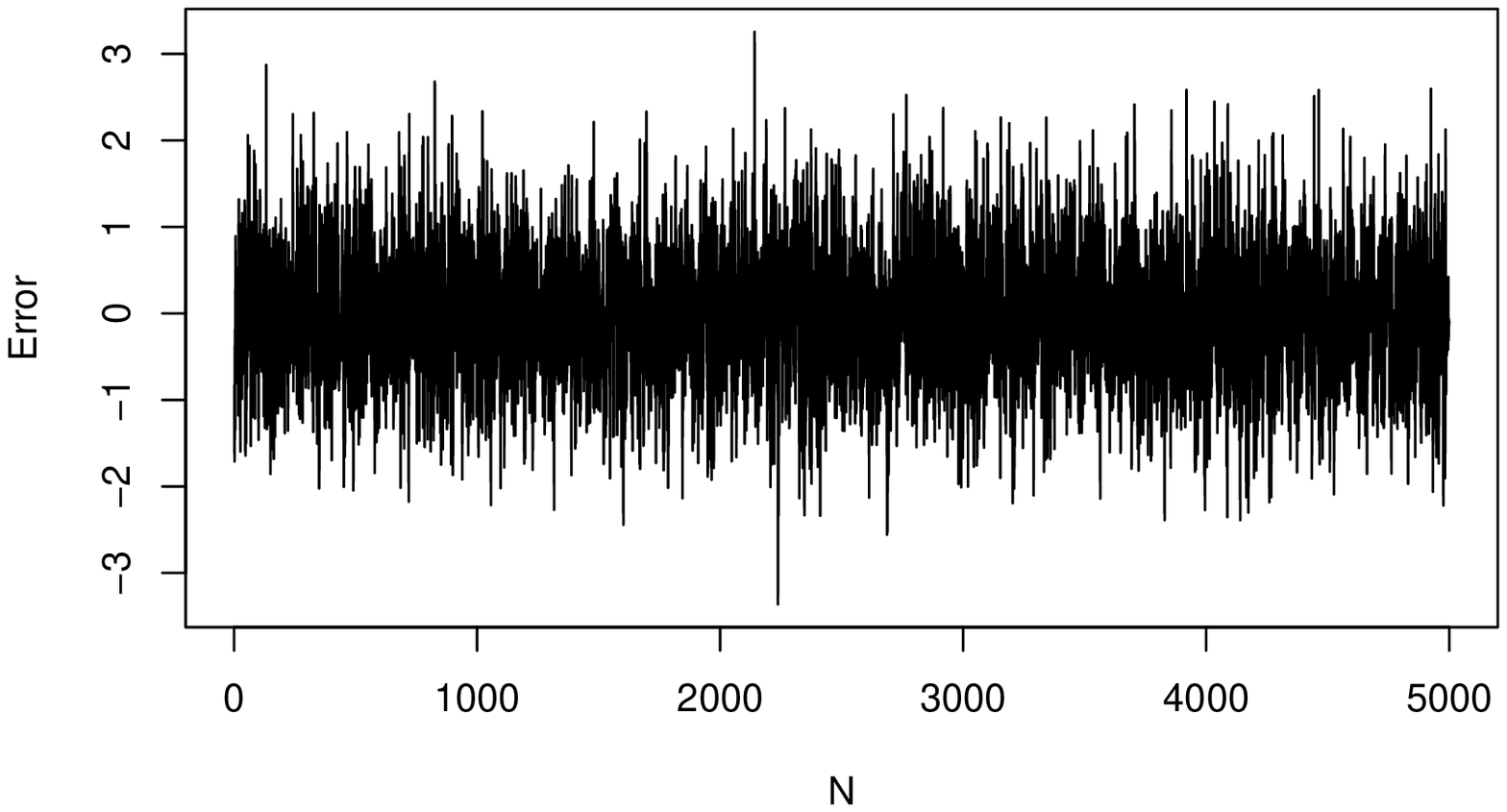}
  \label{erp2}
\end{subfigure}
\caption{Estimation error under different values of $N$ from 1 to 5000 random times. Left: Jittered Sampling, Right: Renewal Sampling}
\label{errorplot}
\end{figure}

Secondly, we plot the histogram of the estimation error defined previously for a fixed value of $N=5000$. 
\begin{figure}[h!]
\centering
\begin{subfigure}{.5\textwidth}
  \centering
  \includegraphics[width=.9\linewidth]{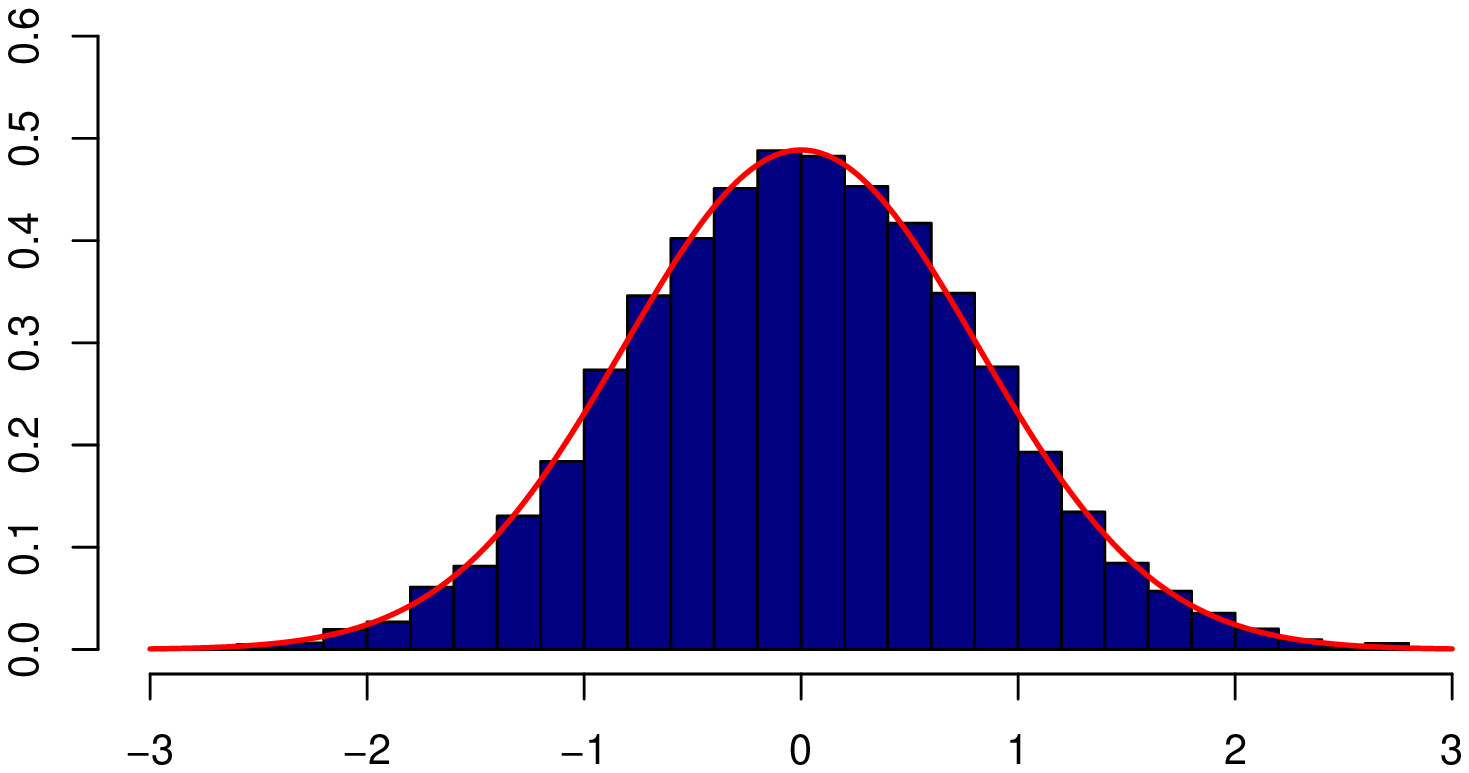}
  \label{ejs_hist}
\end{subfigure}%
\begin{subfigure}{.5\textwidth}
  \centering
  \includegraphics[width=.9\linewidth]{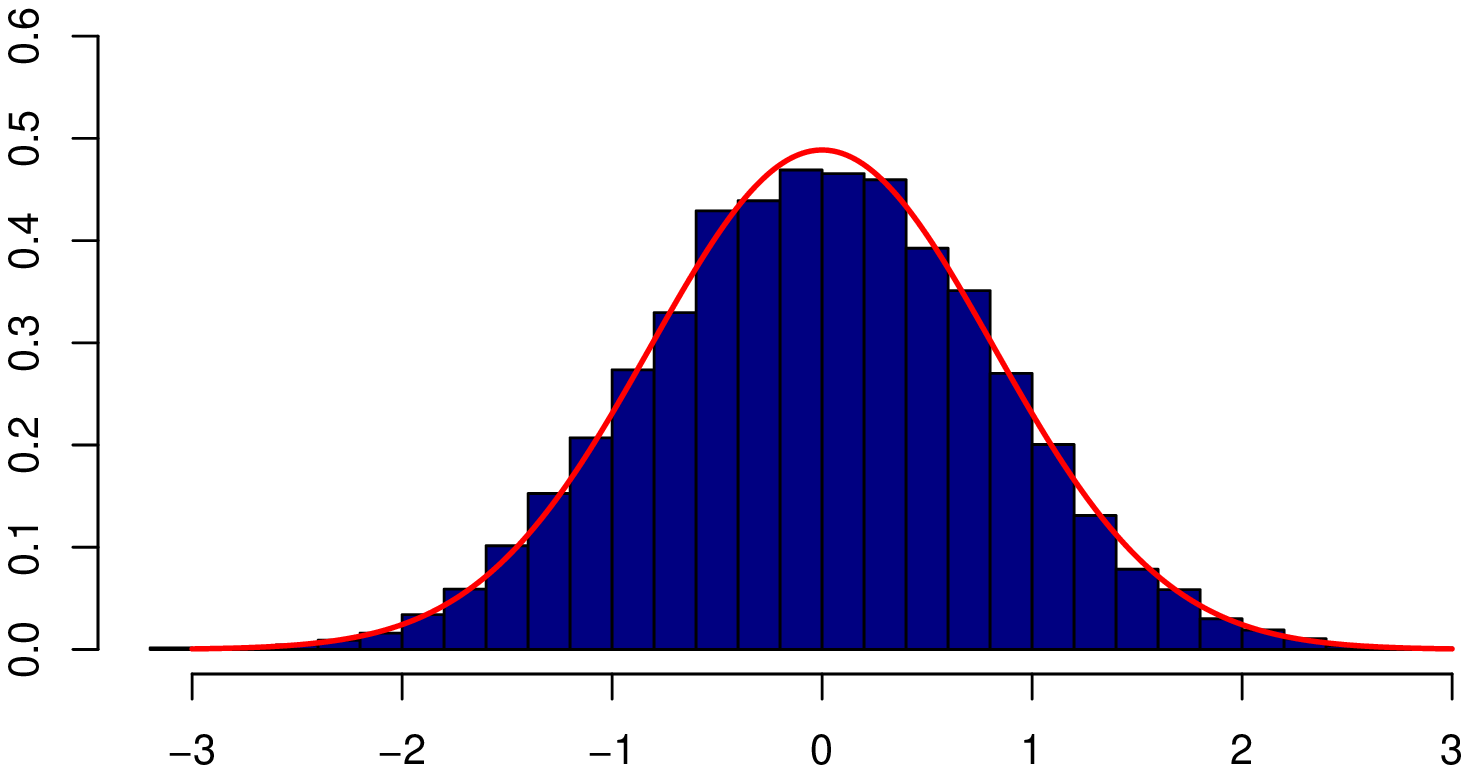}
  \label{erp_hist}
\end{subfigure}
\caption{Histogram of the estimation error under different random times and fixed value of $N=5000$. Left: Jittered Sampling, Right: Renewal Sampling}
\label{errorhist}
\end{figure}

\begin{table}[h!]
\centering
\begin{tabular}{ccc}
\hline \hline 
Random Time & Mean & Variance \\ \hline
Jittered Sampling & 0.003081114 & 0.6664274 \\
Renewal Sampling & -0.000968569 & 0.6666846 \\ 
\hline \hline 
\end{tabular}
\caption{Mean and variance of the estimation error}\label{tab:tab-error}
\end{table}

\pagebreak
\subsection{Conclusions} 
In relation to simulation results previously shown, we can state the following conclusions

\begin{itemize}
\item It is quite easy to check, in Figure \ref{njs1} for jittered sampling and Figure \ref{njs2} for renewal sampling, that as the value of $N$ increases both type of random times are closer to their deterministic or equally spaced version. 
\item When working with renewal observations, we consider an unbounded support so is quite likely that $\tau_{N} > 1$; in order to avoid this, it was necessary to consider $\alpha \approx 1$ such that $\tau_{\alpha N} < 1$. The results given in Table \ref{alpha_tau} exhibits that for large values of $N$ and $\alpha$ close to $1$, it is possible to ensure that $\tau_{\alpha N} < 1$.
\item It is possible to notice that, as the value of $N$ grows, the value of $Q_{N}$ is closer to $1/3$.
\item  In Figure \ref{hist-plot} and in Table \ref{hist-tab}, it is possible to see the asymptotic normality of the renormalized sequence.
\item As we point out in  theorem \ref{an_dn}, the normal distribution with its corresponding parameters is shown in Figure  \ref{hist-plot} and Table \ref{hist-tab}, for both type of random times.
\item  In both Figures \ref{errorplot}, \ref{errorhist} and in Table \ref{tab:tab-error}, for both type of random times, the behavior of the error is as expected, i.e. is centered and with variance not depending on the value of $N$ and besides following a normal distribution. 
\end{itemize}


\section{Appendix}

\subsection{Proofs}

{\bf Proof of Lemma \ref{ll1}:}  Let us separate the proof upon the two situations (\ref{js}) and (\ref{rp}).

\vskip0.2cm

\b{\bf Jittered sampling case: }In this case, as discussed in Section \ref{sec22}, we take $\alpha =1$ and $Y_{\tau_{N}}=0$.  Firstly, we compute the first moment of $A_{N}$, i.e.
\begin{align*}
\mathbb{E} \left[ A_{N} \right] &= \mathbb{E} \left[ \mathbb{E} \left[ A_{N} \mid \tau  \right]  \right] \\
&=  \mathbb{E} \left[ \mathbb{E} \left[ \left.  \dfrac{1}{N} \sum_{j=0}^{N-2} \left( \frac{j+1}{N} + \nu_{j} \right) \left( W_{\frac{j+1}{N} + \nu_{j+1}} - W_{\frac{j}{N} + \nu_{j}} \right) \right|  \nu  \right] \right] \nonumber\\
&= \mathbb{E} \left[ \dfrac{1}{N} \sum_{j=0}^{N-2} \left( \frac{j+1}{N} + \nu_{j+1} \right) \mathbb{E} \left[ W_{\frac{j+1}{N} + \nu_{j+1}} - W_{\frac{j}{N} + \nu_{j}} \mid \tau  \right] \right]\\
&= 0.
\end{align*}
The conditioning with respect to $\tau$ (above and throughout) means that we conditions with respect to the sigma-field generated by the random variables $\nu_{i, N}, i=1,\ldots ,N$. The $L^{2}(\Omega)$ norm of $A_{N}$ can be calculated as follows
\begin{align}
\mathbb{E} \left[ A_{N}^{2} \right] &= Var (A_{N}) \nonumber \\	 
&= \mathbb{E} \left[ Var(A_{N} \mid \tau) \right] + Var(\mathbb{E} \left[ A_{N} \mid \tau  \right]) = \mathbb{E} \left[ Var(A_{N} \mid \tau) \right] \nonumber \\ 
&= \mathbb{E}_{\nu} \left[ \dfrac{1}{N^{2}} \sum_{j=0}^{N-2} \left( \frac{j+1}{N} + \nu_{j+1} \right)^{2} \left( \nu_{j+1} - \nu_{j} + \frac{1}{N} \right)  \right] \nonumber \\ 
&= \mathbb{E} \left[ \dfrac{1}{N^{2}} \sum_{j=0}^{N-2} \left( \frac{(j+1)^{2} \nu_{j+1}}{N^{2}} + \frac{2(j+1) \nu_{j+1}^{2}}{N} + \nu_{j+1}^{3} - \frac{(j+1)^{2} \nu_{j}}{N^{2}} - \frac{2(j+1) \nu_{j} \nu_{j+1}}{N} \right. \right. \nonumber \\ 
&- \left. \left. \nu_{j} \nu_{j+1}^{2} + \frac{(j+1)^{2}}{N^{3}} + \frac{2(j+1) \nu_{j+1}}{N^{2}} + \frac{ \nu_{j+1}^{2}}{N}   \right) \right] \nonumber \\	 
&=: E_{N}^{(1)} + E_{N}^{(2)} + E_{N}^{(3)} - E_{N}^{(4)} - E_{N}^{(5)} - E_{N}^{(6)} + E_{N}^{(7)} + E_{N}^{(8)} + E_{N}^{(9)}. \label{E2-JS}
\end{align}
Notices that   terms $E_{N}^{(1)}, E_{N}^{(4)} , E_{N}^{(5)} , E_{N}^{(6)}$ and $E_{N}^{(8)}$ due to the assumption that $\mathbb{E}(\nu_{i, N})= 0$ for every $i, N$,  $E_{N}^{(3)}$ is also  equal to zero do the symmetry of the law of $\nu _{i, N}$. Therefore
\begin{eqnarray}\label{2d-2}
\mathbb{E} \left[ A_{N}^{2} \right] &=& E_{N}^{(2)}+ E_{N}^{(7)}+  E_{N}^{(9)}.
\end{eqnarray}
We evaluate the three summands in the right-hand side above. By (\ref{2d-1})
\begin{align}
E_{N}^{(2)} &= \mathbb{E} \left[ \dfrac{1}{N^{2}} \sum_{j=0}^{N-2} \dfrac{2 \nu_{j+1}^{2} (j+1)}{N} \right] = \dfrac{2}{N^{3}} \sum_{j=0}^{N-2} (j+1) \mathbb{E} \left[ \nu_{j+1}^{2} \right] \nonumber \\
&= \dfrac{2c_{1}}{ N^{5}}   \left( \dfrac{N(N-1)}{2} \right)=\frac{c_{1}}{N ^ {3}} + o\left( \frac{1}{ N ^ {3}}\right). \label{JS-E2}
\end{align}
Next 
\begin{align}
E_{N}^{(7)} &= \mathbb{E} \left[ \dfrac{1}{N^{2}} \sum_{j=0}^{N-2} \dfrac{(j+1)^{2}}{N^{3}} \right]  = \dfrac{1}{N^{5}} \left( \dfrac{(N-1)(N)(2N-1)}{6} \right) \nonumber \\
&=\frac{1}{3} \frac{1}{N ^ {2}} +o\left(\dfrac{1}{N^{2}}\right) \label{JS-E7}
\end{align}
and, again by (\ref{2d-1}), for $N\geq 2$,
\begin{align}
E_{N}^{(9)} &= \mathbb{E} \left[ \dfrac{1}{N^{2}} \sum_{j=0}^{N-2} \dfrac{\nu_{j+1}^{2}}{N} \right]  = \dfrac{1}{N^{3}} \sum_{j=0}^{N-2} \mathbb{E} \left[ \nu_{j+1}^{2} \right] =  \dfrac{1}{N^{3}} \sum_{j=0}^{N-2}c_{1}  \dfrac{1}{ N^{2}}  = c_{1}\dfrac{1}{(N-1)^{4}}. \label{JS-E9}
\end{align}
From (\ref{2d-2}), (\ref{JS-E2}), (\ref{JS-E7}) and (\ref{JS-E9}), we obtain the conclusion in JS case.

\vskip0.3cm

\b {\bf Renewal sampling case: } Recall that in this case the number of observations is $N_{\alpha}= [\alpha N]$ with $\alpha <1$ close to $1$. As before, by conditioning on $\tau$,
\begin{align*}
\mathbb{E} \left[ A_{N} \right] &= \mathbb{E} \left[ \mathbb{E} \left[ A_{N} \vert \tau \right] \right] = \left. \mathbb{E} \left[ \mathbb{E} \left[ \dfrac{1}{N} \sum_{j=0}^{N_{\alpha}-1} \tau_{j+1} \left( W_{\tau_{j+1}} - W_{\tau_{j}} \right) \right| \tau \right] \right] \\
&=  \mathbb{E} \left[ \left. \dfrac{1}{N} \sum_{j=0}^{N_{\alpha}-1} \tau_{j+1} \mathbb{E} \left[ W_{\tau_{j+1}} - W_{\tau_{j}} \right| \tau \right] \right] = 0,
\end{align*}
and, using that $A_{N}$ has, conditionally on $\tau$, a Gaussian distribution, we obtain
\begin{align}
\mathbb{E} \left[ A_{N}^{2} \right] &= Var (A_{N}) = \mathbb{E} \left[ Var(A_{N} \vert \tau) \right] + Var (\mathbb{E} \left[ A_{N} \vert \tau \right]) \nonumber \\
&= \mathbb{E} \left[ Var(A_{N} \vert \tau) \right] \nonumber = \mathbb{E} \left[ Var \left(  \dfrac{1}{N} \sum_{j=0}^{N_{\alpha}-1} \tau_{j+1} \left( W_{\tau_{j+1}} - W_{\tau_{j}} \right) \right) \right] \nonumber \\
&= \mathbb{E} \left[ \dfrac{1}{N^{2}} \sum_{j=0}^{N_{\alpha}-1} \tau_{j+1}^{2} \left( \tau_{j+1} - \tau_{j} \right) \right]= \mathbb{E} \left[ \dfrac{1}{N^{2}} \sum_{j=0}^{N_{\alpha}-1} \tau_{j+1}^{3} - \dfrac{1}{N^{2}} \sum_{j=0}^{N_{\alpha}-1} \tau_{j+1}^{2} \tau_{j}  \right] \nonumber \\
&= \dfrac{1}{N^{2}} \sum_{j=0}^{N_{\alpha}-1} \mathbb{E} \left[ \tau_{j+1}^{3} \right] - \dfrac{1}{N^{2}} \sum_{j=0}^{N_{\alpha}-1} \mathbb{E} \left[ \tau_{j+1}^{2} \tau_{j}  \right]=: E_{N}^{(1)} - E_{N}^{(2)}. \label{E2-RS}
\end{align}
For the first term of \eqref{E2-RS}, we have
\begin{align}
E_{N}^{(1)} &=  \dfrac{1}{N^{2}} \sum_{j=0}^{N_{\alpha}-1} \mathbb{E} \left[ \tau_{j+1}^{3} \right] = \dfrac{1}{N^{2}} \sum_{j=0}^{N_{\alpha}-1} \int_{0}^{\infty} x^{3} \dfrac{N^{j+1}}{\Gamma(j+1)} x^{j} e^{-Nx} dx\nonumber\\
&= \dfrac{1}{N^{2}} \sum_{j=0}^{N-1} \dfrac{N^{j+1}}{N^{j+4}} \dfrac{\Gamma(j+4)}{\Gamma(j+1)} = \dfrac{1}{N^{2}} \sum_{j=0}^{N_{\alpha}-1}  \dfrac{(j+3)(j+2)(j+1) \Gamma(j+1)}{N^{3} \Gamma(j+1)} \nonumber \\
&= \dfrac{1}{N^{5}} \sum_{j=0}^{N_{\alpha}-1} (j^{3} + 6j^{2} + 11j + 6) \nonumber \\
&= \dfrac{1}{N^{5}} \left( \dfrac{(N_{\alpha}-1)(N_{\alpha})}{2} \right)^{2} + \dfrac{(N_{\alpha}-1)(N_{\alpha})(2(N_{\alpha}-1)+1)}{N^{5}} + \dfrac{11}{2} \dfrac{(N_{\alpha}-1)(N_{\alpha})}{N^{5}} + \dfrac{6N_{\alpha}}{N^{5}} \nonumber \\
&\sim \frac{\alpha^ {4}}{4N}+ \frac{3\alpha^ {3}}{N ^ {2}} +o\left( \frac{1}{N^ {2}}\right).
 \label{RS-E1}
\end{align}

For the second term of \eqref{E2-RS}, we consider the joint density computed by \cite{araya2019}  (see Table (\ref{densities}))
\begin{align}
E_{N}^{(2)} &= \dfrac{1}{N^{2}} \sum_{j=0}^{N_{\alpha}-1} \mathbb{E} \left[ \tau_{j+1}^{2} \tau_{j}  \right]= \dfrac{1}{N^{2}} \sum_{j=0}^{N_{\alpha}-1} \int_{0}^{\infty} \int_{0}^{\tau_{j+1}} \tau_{j} \tau_{j+1}^{2} \dfrac{N^{j+1}}{\Gamma(j)} \tau_{j}^{j-1} e^{-N \tau_{j+1}} d \tau_{j} d \tau_{j+1} \nonumber \\
&= \dfrac{1}{N^{2}} \sum_{j=0}^{N_{\alpha}-1} \dfrac{N^{j+1}}{\Gamma(j)} \int_{0}^{\infty} \tau_{j+1}^{2} e^{-N \tau_{j+1}} \int_{0}^{\tau_{j+1}} \tau_{j}^{j} d \tau_{j} d \tau_{j+1} \nonumber \\
&= \dfrac{1}{N^{2}} \sum_{j=0}^{N_{\alpha}-1} \dfrac{N^{j+1}}{(j+1) \Gamma(j)} \int_{0}^{\infty} \tau_{j+1}^{j+3} e^{-N \tau_{j+1}} d \tau_{j+1}= \dfrac{1}{N^{2}} \sum_{j=0}^{N_{\alpha}-1} \dfrac{N^{j+1}}{(j+1) \Gamma(j)} \dfrac{\Gamma(j+4)}{N^{j+4}} \nonumber \\
&= \dfrac{1}{N^{5}} \sum_{j=0}^{N_{\alpha}-1} \dfrac{(j+3)(j+2)(j+1) j \Gamma(j)}{(j+1) \Gamma(j)}= \dfrac{1}{N^{5}} \sum_{j=0}^{N_{\alpha}-1} (j^{3} + 5j^{2} + 6j) \nonumber \\
&= \dfrac{1}{N^{5}} \left( \dfrac{(N_{\alpha}-1)(N_{\alpha})}{2} \right)^{2} + \dfrac{5}{N^{5}} \left( \dfrac{(N_{\alpha}-1)(N_{\alpha})(2(N_{\alpha}-1)+1)}{6} \right) + \dfrac{6}{N^{5}} \left( \dfrac{(N_{\alpha}-1)(N_{\alpha})}{2} \right) \nonumber \\
&\sim \dfrac{\alpha ^ {4}}{4N} + \dfrac{7\alpha ^ {3}}{6N^{2}} +o\left( \frac{1}{N ^ {2}}\right). \label{RS-E2}
\end{align}
Replacing (\ref{RS-E1}) and (\ref{RS-E2}) in (\ref{E2-RS})
\begin{align}
\mathbb{E} \left[ A_{N}^{2} \right] &\sim  \dfrac{\alpha ^ {3}}{3N^{2}} + o\left( \frac{1}{N^ {2}}\right). \label{cte-rs}
\end{align}
\qed

\vskip0.2cm

{\bf Proof of Proposition \ref{L2-conv}: }Again, we separately discuss the two cases of random sampling.

\b{\bf Jittered sampling case: } Recall that in this case $\alpha =1$.  First, we compute the first moment of \eqref{qn}, i.e
\begin{align*}
\mathbb{E} \left[ Q_{N} \right] &= \mathbb{E} \left[ \sum_{j=0}^{N-1} \left( \frac{j+1}{N} + \frac{X_{j+1}}{N} \right)^{2} \left( \frac{X_{j+1}}{N} - \frac{X_{j}}{N} + \frac{1}{N}  \right)   \right] \\
&= \dfrac{1}{N^{3}} \sum_{j=0}^{N-1} \mathbb{E} \left[ (j+1)^{2} X_{j+1} - (j+1)^{2}X_{j} + (j+1)^{2} + 2(j+1)X_{j+1}^{2} - 2(j+1)X_{j} X_{j+1} \right. \\
& \left.+  2(j+1)X_{j+1} + X_{j+1}^{3} - X_{j}X_{j+1}^{2} + X_{j+1}^{2}  \right] \\
&=  \dfrac{1}{N^{3}} \sum_{j=0}^{N-1} \left[ Q_{N}^{(1)} - Q_{N}^{(2)} + Q_{N}^{(3)} + Q_{N}^{(4)} - Q_{N}^{(5)} + Q_{N}^{(6)} + Q_{N}^{(7)} - Q_{N}^{(8)} + Q_{N}^{(9)} \right].
\end{align*}
From (\ref{3d-1}) and the hypothesis (\ref{2d-1}) on $\tau$, we have
\begin{equation*}
\mathbb{E}( X_{j})= 0 \mbox{ and } \mathbb{E} (X_{j}^ {2}) = c_{1}.
\end{equation*}
Thus $Q_{N}^{(1)}$, $Q_{N}^{(2)}$, $Q_{N}^{(5)}$ , $Q_{N}^{(6)}$, $Q_{N}^{(7)}$ and $Q_{N}^{(8)}$ are equal to zero, while the remaining terms can be computed as follows
\begin{align}
Q_{N}^{(3)} &=  \mathbb{E} \left[ \dfrac{1}{N^{3}} \sum_{j=0}^{N-1} (j+1)^{2} \right] = \dfrac{N(N+1)(2N+1)}{6} \nonumber \\
&= \dfrac{2N^{3} + 3N^{2} + N}{6} \xrightarrow[N \to \infty]{}  \dfrac{1}{3} \label{q3}
\end{align}

\begin{align}
Q_{N}^{(4)} &=  \mathbb{E} \left[ \dfrac{1}{N^{3}} \sum_{j=0}^{N-1} 2(j+1)X_{j+1}^{2} \right] = \dfrac{2}{N^{3}} \sum_{j=0}^{N-1} (j+1) Var \left( X_{j+1}^{2} \right) \nonumber \\ 
&= \dfrac{2}{N^{3}} \sum_{j=0}^{N-1} (j+1) \dfrac{1}{12}= \dfrac{1}{6 N^{3}}  \left( \dfrac{N(N+1)}{2} \right)  \xrightarrow[N \to \infty]{} 0 \label{q4}
\end{align}
and 
\begin{align}
Q_{N}^{(9)} &=  \mathbb{E} \left[ \dfrac{1}{N^{3}} \sum_{j=0}^{N-1} X_{j+1}^{2} \right]= \dfrac{1}{N^{3}} \sum_{j=0}^{N-1} Var \left( X_{j+1} \right) \nonumber \\
&= \dfrac{1}{N^{3}} \sum_{j=0}^{N-1} c_{1}\xrightarrow[N \to \infty]{} 0 \label{q9}
\end{align}
Taking into account \eqref{q3}, \eqref{q4} and \eqref{q9} we conclude
\begin{equation} \label{lim-qj}
\lim_{N \to \infty} \mathbb{E} \left[ Q_{N} \right] \xrightarrow[N \to \infty]{} \dfrac{1}{3}.
\end{equation}
Secondly, we study the second moment of  \eqref{qn}, i.e.
\begin{align*}
\mathbb{E} \left[ Q_{N}^{2} \right] &= \dfrac{1}{N^{6}} \sum_{j=0}^{N-1} \mathbb{E} \left[ (j+1 + X_{j+1})^{4} (X_{j+1} - X_{j} + 1)^{2} \right] \\
&+ \dfrac{1}{N^{6}} \sum_{0 \leq j \neq k \leq N-1} \mathbb{E} \left[ (j+1 + X_{j+1})^{2} (k+1 + X_{k+1})^{2} (X_{j+1} - X_{j} + 1) (X_{k+1} - X_{k} + 1) \right] \\
&:= Q_{N}^{(I)} + Q_{N}^{(II)}.
\end{align*}
Some of the summands that compose $Q_{N}^ {(I)}$ are zero (those involving odd order moments of $X_{j}$ or $X_{j+1}$, which vanish) and the other summands of $Q_{N} ^ {(I)} $ converge to zero. For example
\begin{align*}
\mathbb{E} \left[ \dfrac{1}{N^{6}} \sum_{j=0}^{N-1}  (j+1)^{4} X_{j+1}^{2} \right] &= \dfrac{1}{N^{6}} \sum_{j=0}^{N-1} (j+1)^{4} Var(X_{j+1}) \\
&= \dfrac{c_{1}}{ N^{6}} \left( \dfrac{6N^{5} + 15N^{4} + 10N^{3} - N}{30} \right) \xrightarrow[N \to \infty]{} 0
\end{align*}
or
\begin{align*}
\mathbb{E} \left[ \dfrac{1}{N^{6}} \sum_{j=0}^{N-1} 2(j+1)^{4} X_{j} X_{j+1}  \right] &= \dfrac{2}{N^{6}} \sum_{j=0}^{N-1} (j+1) \mathbb{E} \left[ X_{j} X_{j+1} \right] =0
\end{align*}
while the remaining terms can be computed in the same way as the last two. \\

The summand $ Q_{N}^{(II)}$ give the limit of $Q_{N}$. Actually, the only non-vanishing term in  $ Q_{N}^{(II)}$ is the one not depending on $X_{j}$, i.e.
\begin{align*}
\mathbb{E} \left[ \dfrac{1}{N^{6}} \sum_{0 \leq j \neq k \leq N-1} (j+1)^{2} (k+1)^{2} \right] &= \dfrac{1}{N^{6}} \sum_{0 \leq j \neq k \leq N-1} (j+1)^{2} (k+1)^{2}.
\end{align*}
To compute this term, we use the following identity
\begin{align}
\sum_{1 \leq i \neq j \leq N} a_{i} a_{j}  &= \left( \sum_{i=1}^{N} a_{i} \right)^{2}  - \sum_{i=1}^{N} a_{i}^{2} \label{suma}
\end{align}

Therefore,
\begin{align}
 &\dfrac{1}{N^{6}} \sum_{0 \leq j \neq k \leq N-1} (j+1)^{2} (k+1)^{2}  = \dfrac{1}{N^{6}} \left[ \left( \sum_{j=0}^{N-1} (j+1)^{2} \right)^{2} \sum_{j=0}^{N-1} (j+1)^{4} \right] \nonumber \\
 &= \dfrac{1}{N^{6}} \left[ \left( \dfrac{N(N+1)(2N+1)}{6} \right)^{2} - \left( \dfrac{6N^{5} + 15N^{4} + 10N^{3} - N}{30} \right) \right] \nonumber \\
 &= \dfrac{1}{N^{6}} \left[ \dfrac{(2N^{3} + 3N^{2} + N)^{2}}{36} - \dfrac{6N^{5} - 15N^{4} - 10N^{3} + N}{30} \right] \xrightarrow[N \to \infty]{}  \dfrac{1}{9}. \label{e2} 
\end{align}
Taking into account \eqref{lim-qj} and \eqref{e2}, we obtain the conclusion.

\vskip0.3cm

\b {\bf Renewal sampling case: } As before, we start by computing the expectation of $Q_{N}$. We have
\begin{align}
\mathbb{E} \left[ Q_{N} \right] &= \mathbb{E} \left[ \sum_{j=0}^{N_{\alpha}-1} \tau_{j+1}^{2} \left( \tau_{j+1} - \tau_{j} \right) \right] \nonumber \\
&= \mathbb{E} \left[ \sum_{j=0}^{N_{\alpha}-1} \tau_{j+1}^{3} \right] - \mathbb{E} \left[\sum_{j=0}^{N_{\alpha}-1} \tau_{j+1}^{2} \tau_{j}  \right]= Q_{N}^{(1)} - Q_{N}^{(2)} \label{eqn}
\end{align}
with
\begin{align}
Q_{N}^{(1)} &= \mathbb{E} \left[ \sum_{j=0}^{N-1} \tau_{j+1}^{3} \right]= \sum_{j=0}^{N_{\alpha}-1} \int_{0}^{\infty} \tau_{j+1}^{3} \dfrac{N^{j+1}}{\Gamma(j+1)} \tau_{j+1}^{j} e^{-N \tau_{j+1}} d \tau_{j+1} \nonumber \\
&= \sum_{j=0}^{N_{\alpha}-1} \dfrac{N^{j+1}}{\Gamma(j+1)} \dfrac{\Gamma(j+4)}{N^{j+4}}= \dfrac{1}{N^{3}} \sum_{j=0}^{N_{\alpha}-1} (j+3)(j+2)(j+1) \nonumber \\
&= \dfrac{1}{N^{3}} \sum_{j=0}^{N_{\alpha}-1} [j^{3} + 6j^{2} + 11j + 6] \label{eqn1}
\end{align}
and

\begin{align}
Q_{N}^{(2)} &= \mathbb{E} \left[ \sum_{j=0}^{N_{\alpha}-1} \tau_{j+1}^{2} \tau_{j}  \right] = \sum_{j=0}^{N_{\alpha}-1} \int_{0}^{\infty} \int_{0}^{\tau_{j+1}} \tau_{j+1}^{2} \tau_{j} \dfrac{N^{j+1}}{\Gamma(j)} \tau_{j}^{j-1} e^{-N \tau_{j+1}} d \tau_{j} d \tau_{j+1} \nonumber \\
&= \sum_{j=0}^{N_{\alpha}-1} \dfrac{N^{j+1}}{\Gamma(j)} \int_{0}^{\infty} \tau_{j+1}^{2} e^{-N \tau_{j+1}}  \int_{0}^{\tau_{j+1}} \tau_{j}^{j} d \tau_{j} d \tau_{j+1}= \sum_{j=0}^{N_{\alpha}-1} \dfrac{N^{j+1}}{\Gamma(j) (j+1)} \int_{0}^{\infty} \tau_{j+1}^{j+4-1} e^{-N \tau_{j+1}} d \tau_{j+1} \nonumber \\
&= \sum_{j=0}^{N_{\alpha}-1}  \dfrac{N^{j+1}}{N^{j+4}} \dfrac{\Gamma(j+4)}{\Gamma(j) (j+1)} = \dfrac{1}{N^{3}} \sum_{j=0}^{N_{\alpha}-1} (j+3)(j+2)j= \dfrac{1}{N^{3}} \sum_{j=0}^{N_{\alpha}-1} [j^{3} - 5j^{2} +6j]  \label{eqn2}.
\end{align}

Replacing \eqref{eqn1} and \eqref{eqn2} in \eqref{eqn}, it results
\begin{align*}
\mathbb{E} \left[ Q_{N} \right] &= \dfrac{1}{N^{3}} \sum_{j=0}^{N_{\alpha}-1} [j^{2} - 5j +6] 
= \dfrac{1}{N^{3}}  \left[ \dfrac{2N_{\alpha}^{3} - 18N_{\alpha}^{2} + 52N_{\alpha}}{6}  \right]  \xrightarrow[N \to \infty]{}  \dfrac{\alpha ^ {3}}{3}
\end{align*}

For the second moment of  $Q_{N}$ we have
\begin{align}
\mathbb{E} \left[ Q_{N}^{2} \right] &= \mathbb{E} \left[ \left( \sum_{j=0}^{N_{\alpha}-1} \tau_{j+1}^{2} \left( \tau_{j+1} - \tau_{j} \right) \right)^{2} \right] \nonumber \\
&= \mathbb{E} \left[ \sum_{j=0}^{N_{\alpha}-1} \tau_{j+1}^{4} \left( \tau_{j+1} - \tau_{j}  \right)^{2} \right] + 2 \mathbb{E} \left[ \sum_{j<k}^{N_{\alpha}-1} \tau_{j+1}^{2} \tau_{k+1}^{2} \left( \tau_{j+1} - \tau_{j} \right) \left( \tau_{k+1} - \tau_{k} \right)  \right] \nonumber \\ 
&= \mathbb{E} \left[ \sum_{j=0}^{N_{\alpha}-1} \tau_{j+1}^{6} - 2 \sum_{j=0}^{N_{\alpha}-1} \tau_{j+1}^{5} \tau_{j} + \sum_{j=0}^{N_{\alpha}-1} \tau_{j+1}^{4} \tau_{j}^{2}  \right] \nonumber \\
&+ 2 \mathbb{E} \left[ \sum_{j<k}^{N_{\alpha}-1} \tau_{j+1}^{3} \tau_{k+1}^{3} - \sum_{j<k}^{N_{\alpha}-1} \tau_{j+1}^{3} \tau_{k} \tau_{k+1}^{2} - \sum_{j<k}^{N_{\alpha}-1} \tau_{j+1}^{2} \tau_{j} \tau_{k+1}^{3} + \sum_{j<k}^{N_{\alpha}-1} \tau_{j} \tau_{j+1}^{2} \tau_{k} \tau_{k+1}^{2} \right] \nonumber \\
&:= Q_{N}^{(I)} - Q_{N}^{(II)} + Q_{N}^{(III)} + Q_{N}^{(IV)} - Q_{N}^{(V)} - Q_{N}^{(VI)} + Q_{N}^{(VII)} \label{eqnn}
\end{align}
and the above terms can be calculated by using the joint densities from Table \ref{densities}.  We calculate first the sum $Q_{N}^{(I)} - Q_{N}^{(II)} + Q_{N}^{(III)} $.  First
\begin{align}
Q_{N}^{(I)} &= \mathbb{E} \left[ \sum_{j=0}^{N_{\alpha}-1} \tau_{j+1}^{6} \right]=\sum_{j=0}^{N_{\alpha}-1}  \int_{0}^{\infty} \tau_{j+1}^{6} \dfrac{N^{j+1}}{\Gamma(j+1)} \tau_{j+1}^{j} e^{-N \tau_{j+1}} d \tau_{j+1} \nonumber \\
&= \sum_{j=0}^{N_{\alpha}-1} \dfrac{N^{j+1}}{\Gamma(j+1)} \int_{0}^{\infty} \tau_{j+1}^{j+7-1} e^{-N \tau_{j+1}} d \tau_{j+1}= \sum_{j=0}^{N_{\alpha}-1} \dfrac{N^{j+1}}{N^{j+7}} \dfrac{\Gamma(j+7)}{\Gamma(j+1)} \nonumber \\
&= \dfrac{1}{N^{6}} \sum_{j=0}^{N_{\alpha}-1} (j+6)(j+5)(j+4)(j+3)(j+2)(j+1) \nonumber \\
&= \dfrac{1}{N^{6}} \sum_{j=0}^{N_{\alpha}-1} [j^{6} + 21j^{5} + 175j^{4} + 735j^{3} + 1624j^{2} + 1764j + 720]. \label{e2-1}
\end{align}
For $Q_{N}^{(II)}, Q_{N}^{(III)}$ we use the joint density of the random vector $(\tau_{j}, \tau_{k})$ 
\begin{align}
Q_{N}^{(II)} &= 2 \mathbb{E} \left[ \sum_{j=0}^{N_{\alpha}-1} \tau_{j+1}^{5} \tau_{j} \right]= 2 \sum_{j=0}^{N_{\alpha}-1} \int_{0}^{\infty} \int_{0}^{\tau_{j+1}} \tau_{j+1}^{5} \tau_{j} \dfrac{N^{j+1}}{\Gamma(j)} \tau_{j}^{j-1} e^{-N \tau_{j+1}} d \tau_{j}d \tau_{j+1} \nonumber \\
&= 2 \sum_{j=0}^{N_{\alpha}-1} \dfrac{N^{j+1}}{\Gamma(j)} \int_{0}^{\infty}  \tau_{j+1}^{5} e^{-N \tau_{j+1}} \int_{0}^{\tau_{j+1}} \tau_{j}^{j} d \tau_{j} d \tau_{j+1} \nonumber \\
&= 2 \sum_{j=0}^{N_{\alpha}-1} \dfrac{N^{j+1}}{\Gamma(j) (j+1)} \int_{0}^{\infty}  \tau_{j+1}^{j+7-1} e^{-N \tau_{j+1}} d \tau_{j+1} \nonumber \\
&= 2 \sum_{j=0}^{N_{\alpha}-1} \dfrac{N^{j+1}}{N^{j+7}} \dfrac{\Gamma(j+7)}{\Gamma(j) (j+1)} = \dfrac{2}{N^{6}} \sum_{j=0}^{N_{\alpha}-1} (j+6)(j+5)(j+4)(j+3)(j+2)j \nonumber \\
&= \dfrac{1}{N^{6}} \sum_{j=0}^{N_{\alpha}-1}  [2j^{6} + 40j^{5} + 310j^{4} + 1160j^{3} + 2088j^{2} + 1440j] \label{e2-2}
\end{align}
and

\begin{align}
Q_{N}^{(III)} &= \mathbb{E} \left[ \sum_{j=0}^{N_{\alpha}-1} \tau_{j+1}^{4} \tau_{j}^{2}  \right] = \sum_{j=0}^{N_{\alpha}-1}  \int_{0}^{\infty} \int_{0}^{\tau_{j+1}} \tau_{j+1}^{4} \tau_{j}^{2} \dfrac{N^{j+1}}{\Gamma(j)} \tau_{j}^{j-1} e^{-N \tau_{j+1}} d \tau_{j}d \tau_{j+1} \nonumber \\
&= \sum_{j=0}^{N_{\alpha}-1} \dfrac{N^{j+1}}{\Gamma(j)} \int_{0}^{\infty}  \tau_{j+1}^{4} e^{-N \tau_{j+1}} \int_{0}^{\tau_{j+1}} \tau_{j}^{j+1} d \tau_{j} d \tau_{j+1} \nonumber \\
&= \sum_{j=0}^{N_{\alpha}-1} \dfrac{N^{j+1}}{\Gamma(j) (j+2)} \int_{0}^{\infty}  \tau_{j+1}^{j+7-1} e^{-N \tau_{j+1}} d \tau_{j+1} \nonumber \\
&= \sum_{j=0}^{N_{\alpha}-1} \dfrac{N^{j+1}}{N^{j+7}} \dfrac{\Gamma(j+7)}{\Gamma(j) (j+2)}  = \dfrac{1}{N^{6}} \sum_{j=0}^{N_{\alpha}-1} (j+6)(j+5)(j+4)(j+3)(j+1)j \nonumber \\
&= \dfrac{1}{N^{6}} \sum_{j=0}^{N_{\alpha}-1} [j^{6} + 19j^{5} + 137j^{4} + 461j^{3} + 702j^{2} + 360j]. \label{e2-3}
\end{align}
By putting together \eqref{e2-1}, \eqref{e2-2} and \eqref{e2-3},
\begin{align}
Q_{N}^{(I)} - Q_{N}^{(II)} + Q_{N}^{(III)} = \dfrac{1}{N^{6}} \sum_{j=0}^{N_{\alpha}-1} [j^{4} + 18j^{3} + 119j^{2} + 342j + 360] \xrightarrow[N \to \infty]{} 0 \label{lim1},
\end{align} 
 For terms $Q_{N}^{(IV)}, Q_{N}^{(V)}, Q_{N}^{(VI)}$ and $Q_{N}^{(VII)}$ it is necessary to consider the following joint densities $f_{\tau_{j+1} , \tau_{k+1}}$, $f_{\tau_{j}, \tau_{k} , \tau_{l}}$ and $f_{\tau_{j}, \tau_{k} , \tau_{l}, \tau_{m}}$ with all the different indices. We also use the identity
\begin{equation}
\int_{0}^{b} x^{a} (b-x)^{c} dx = \dfrac{\Gamma(a+1) \Gamma(c+1)}{\Gamma(a+c+2)} b^{a+c+1} \label{inte}.
\end{equation}
We have the following calculations
\begin{align}
Q_{N}^{(IV)} &= \mathbb{E} \left[ \sum_{j<k}^{N_{\alpha}-1} \tau_{j+1}^{3} \tau_{k+1}^{3} \right] \nonumber \\
&= \sum_{j<k}^{N_{\alpha}-1} \int_{0}^{\infty} \int_{0}^{\tau_{k+1}} \tau_{j+1}^{3} \tau_{k+1}^{3} \dfrac{N^{k+1}}{\Gamma(j+1) \Gamma(k-j)} \tau_{j+1}^{j} (\tau_{k+1} - \tau_{j+1})^{k-j-1} e^{-N \tau_{k+1}} d \tau_{j+1} d \tau_{k+1} \nonumber \\
&= \sum_{j<k}^{N_{\alpha}-1} \dfrac{N^{k+1}}{\Gamma(j+1) \Gamma(k-j)} \int_{0}^{\infty} \tau_{k+1}^{3} e^{-N \tau_{k+1}} \int_{0}^{\tau_{k+1}} \tau_{j+1}^{j+3} (\tau_{k+1} - \tau_{j+1})^{k-j-1}   d \tau_{j+1} d \tau_{k+1} \nonumber \\
&= \sum_{j<k}^{N_{\alpha}-1} \dfrac{N^{k+1}}{\Gamma(j+1) \Gamma(k-j)} \dfrac{\Gamma(j+4) \Gamma(k-j)}{\Gamma(k+4)} \int_{0}^{\infty} \tau_{k+1}^{k+7-1} e^{-N \tau_{k+1}} d \tau_{k+1} \nonumber \\
&= \sum_{j<k}^{N_{\alpha}-1} \dfrac{N^{k+1}}{N^{k+7}} \dfrac{\Gamma(j+4) \Gamma(k+7)}{\Gamma(j+1) \Gamma(k+4)}= \dfrac{1}{N^{6}} \sum_{j<k}^{N_{\alpha}-1} (j+3)(j+2)(j+1) (k+6)(k+5)(k+4) \label{e2-4},
\end{align}
and, via \eqref{inte}
\begin{align}
Q_{N}^{(V)} &= \mathbb{E} \left[ \sum_{j<k}^{N_{\alpha}-1} \tau_{j+1}^{3} \tau_{k} \tau_{k+1}^{2} \right] \nonumber \\
&= \sum_{j<k}^{N_{\alpha}-1} \int_{0}^{\infty} \int_{0}^{\tau_{k+1}} \int_{0}^{\tau_{k}} \tau_{j+1}^{3} \tau_{k} \tau_{k+1}^{2} \dfrac{N^{k+1}}{\Gamma(j+1) \Gamma(k-j-1)} \tau_{j+1}^{j} (\tau_{k} - \tau_{j+1})^{k-j-2} e^{-N \tau_{k+1}} d \tau_{j+1} d \tau_{k} d \tau_{k+1} \nonumber \\
&= \sum_{j<k}^{N_{\alpha}-1} \dfrac{N^{k+1}}{\Gamma(j+1) \Gamma(k-j-1)}  \int_{0}^{\infty} \tau_{k+1}^{2} e^{-N \tau_{k+1}} \int_{0}^{\tau_{k+1}} \tau_{k} \int_{0}^{\tau_{k}} \tau_{j+1}^{j+3} (\tau_{k} - \tau_{j+1})^{k-j-2}  d \tau_{j+1} d \tau_{k} d \tau_{k+1} \nonumber \\
&= \sum_{j<k}^{N_{\alpha}-1} \dfrac{N^{k+1}}{\Gamma(j+1) \Gamma(k-j-1)} \dfrac{\Gamma(j+4) \Gamma(k-j-1)}{\Gamma(k+3)}  \int_{0}^{\infty} \tau_{k+1}^{2} e^{-N \tau_{k+1}} \int_{0}^{\tau_{k+1}} \tau_{k}^{k+3} d \tau_{k} d \tau_{k+1} \nonumber \\
&=  \sum_{j<k}^{N_{\alpha}-1} \dfrac{N^{k+1} \Gamma(j+4}{\Gamma(j+1) \Gamma(k+3) (k+4)} \int_{0}^{\infty} \tau_{k+1}^{k+7-1} e^{-N \tau_{k+1}} d \tau_{k+1} \nonumber \\
&= \sum_{j<k}^{N_{\alpha}-1} \dfrac{N^{k+1}}{N^{k+7}} \dfrac{\Gamma(j+4) \Gamma(k+7)}{\Gamma(j+1) \Gamma(k+3) (k+4)}= \dfrac{1}{N^{6}} \sum_{j<k}^{N_{\alpha}-1}  (j+3)(j+2)(j+1)(k+6)(k+5)(k+3) \label{e2-5}
\end{align}
and
\begin{align}
Q_{N}^{(VI)} &= \mathbb{E} \left[ \sum_{j<k}^{N_{\alpha}-1} \tau_{j+1}^{2} \tau_{j} \tau_{k+1}^{3} \right] \nonumber \\
&= \sum_{j<k}^{N_{\alpha}-1} \int_{0}^{\infty} \int_{0}^{\tau_{k+1}} \int_{0}^{\tau_{j+1}} \tau_{j+1}^{2} \tau_{j} \tau_{k+1}^{3} \dfrac{N^{k+1}}{\Gamma(j) \Gamma(k-j)} \tau_{j}^{j-1} (\tau_{k+1} - \tau_{j+1})^{k-j-1} e^{-N \tau_{k+1}} d \tau_{j} d \tau_{j+1} d \tau_{k+1} \nonumber \\
&= \sum_{j<k}^{N_{\alpha}-1} \dfrac{N^{k+1}}{\Gamma(j) \Gamma(k-j)} \int_{0}^{\infty} \tau_{k+1}^{3} e^{-N \tau_{k+1}}   \int_{0}^{\tau_{k+1}} \tau_{j+1}^{2} (\tau_{k+1} - \tau_{j+1})^{k-j-1}  \int_{0}^{\tau_{j+1}} \tau_{j}^{j}  d \tau_{j} d \tau_{j+1} d \tau_{k+1} \nonumber \\
&= \sum_{j<k}^{N_{\alpha}-1} \dfrac{N^{k+1}}{\Gamma(j) \Gamma(k-j) (j+1)} \int_{0}^{\infty} \tau_{k+1}^{3} e^{-N \tau_{k+1}}   \int_{0}^{\tau_{k+1}} \tau_{j+1}^{j+3} (\tau_{k+1} - \tau_{j+1})^{k-j-1} d \tau_{j+1} d \tau_{k+1} \nonumber \\
&= \sum_{j<k}^{N_{\alpha}-1} \dfrac{N^{k+1}}{\Gamma(j) \Gamma(k-j) (j+1)} \dfrac{\Gamma(j+4) \Gamma(k-j)}{\Gamma(k+4)} \int_{0}^{\infty} \tau_{k+1}^{k+7-1} e^{-N \tau_{k+1}} d \tau_{k+1} \nonumber \\
&= \sum_{j<k}^{N_{\alpha}-1} \dfrac{N^{k+1}}{N^{k+7}} \dfrac{\Gamma(j+4) \Gamma(k+7)}{\Gamma(j) (j+1) \Gamma(k+4)} = \dfrac{1}{N^{6}} \sum_{j<k}^{N_{\alpha}-1} (j+3)(j+2)j (k+6)(k+5)(k+4) \label{e2-6}
\end{align}
and finally
\begin{align}
Q_{N}^{(VII)} &= \mathbb{E} \left[ \sum_{j<k}^{N_{\alpha}-1} \tau_{j} \tau_{j+1}^{2} \tau_{k} \tau_{k+1}^{2} \right] \nonumber \\
&= \sum_{j<k}^{N_{\alpha}-1} \int_{0}^{\infty} \int_{0}^{\tau_{k+1}} \int_{0}^{\tau_{k}} \int_{0}^{\tau_{j+1}} \tau_{j} \tau_{j+1}^{2} \tau_{k} \tau_{k+1}^{2} \dfrac{N^{k+1}}{\Gamma(j) \Gamma(k-j-1)} \tau_{j}^{j-1} (\tau_{k} - \tau_{j+1})^{k-j-2} e^{-N \tau_{k+1}} d \tau_{j}d \tau_{j+1} d \tau_{k} d \tau_{k+1} \nonumber \\
&= \sum_{j<k}^{N_{\alpha}-1}  \dfrac{N^{k+1}}{\Gamma(j) \Gamma(k-j-1)}  \int_{0}^{\infty} \tau_{k+1}^{2} e^{-N \tau_{k+1}}  \int_{0}^{\tau_{k+1}}  \tau_{k}     \int_{0}^{\tau_{k}} \tau_{j+1}^{2}   (\tau_{k} - \tau_{j+1})^{k-j-2}  \int_{0}^{\tau_{j+1}}  \tau_{j}^{j}  d \tau_{j} d \tau_{j+1} d \tau_{k} d \tau_{k+1} \nonumber \\
&= \sum_{j<k}^{N_{\alpha}-1}  \dfrac{N^{k+1}}{\Gamma(j) \Gamma(k-j-1) (j+1)}  \int_{0}^{\infty} \tau_{k+1}^{2} e^{-N \tau_{k+1}}  \int_{0}^{\tau_{k+1}}  \tau_{k}     \int_{0}^{\tau_{k}} \tau_{j+1}^{j+3}   (\tau_{k} - \tau_{j+1})^{k-j-2} d \tau_{j+1} d \tau_{k} d \tau_{k+1} \nonumber \\
&= \sum_{j<k}^{N_{\alpha}-1}  \dfrac{N^{k+1}}{\Gamma(j) \Gamma(k-j-1) (j+1)} \dfrac{\Gamma(j+4) \Gamma(k-j-1)}{\Gamma(k+3} \int_{0}^{\infty} \tau_{k+1}^{2} e^{-N \tau_{k+1}}  \int_{0}^{\tau_{k+1}}  \tau_{k+3}   d \tau_{k} d \tau_{k+1} \nonumber \\
&= \sum_{j<k}^{N_{\alpha}-1}  \dfrac{N^{k+1} \Gamma(j+4)}{\Gamma(j) (j+1) \Gamma(k+3) (k+4)}  \int_{0}^{\infty} \tau_{k+1}^{k-7-1} e^{-N \tau_{k+1}} d \tau_{k+1} \nonumber \\
&= \sum_{j<k}^{N_{\alpha}-1} \dfrac{N^{k+1}}{N^{k+7}}  \dfrac{\Gamma(j+4) \Gamma(k+7)}{\Gamma(j) (j+1) \Gamma(k+3) (k+4)} = \dfrac{1}{N^{6}} \sum_{j<k}^{N_{\alpha}-1} (j+3)(j+2)j (k+6)(k+5)(k+3) \label{e2-7}. 
\end{align}
Taking into account the results from \eqref{e2-4}, \eqref{e2-5}, \eqref{e2-6} and \eqref{e2-7}, we get
\begin{align}
Q_{N}^{(IV)} - Q_{N}^{(V)} - Q_{N}^{(VI)} + Q_{N}^{(VII)} &= \dfrac{2}{N^{6}} \sum_{j<k}^{N-1} \left[ (j+3)(j+2)(j+1) (k+6)(k+5)(k+4) \right. \nonumber \\
& - (j+3)(j+2)(j+1)(k+6)(k+5)(k+3) \nonumber \\
& - (j+3)(j+2)j (k+6)(k+5)(k+4)  \nonumber \\
&+ \left. (j+3)(j+2)j (k+6)(k+5)(k+3) \right] \nonumber \\
&= \dfrac{2}{N^{6}} \sum_{j<k}^{N-1} \left[  (j+3)(j+2)(j+1) (k+6)(k+5) \left[ (k+4) - (k+3) \right] \right. \nonumber \\
&+ \left. (j+3)(j+2)j (k+6)(k+5) \left[ (k+3) - (k+4) \right] \right] \nonumber \\
&= \dfrac{2}{N^{6}} \sum_{j<k}^{N-1} \left[ (j+3)(j+2) (k+6)(k+5) \left[ (j+1) - j \right] \right] \nonumber \\
&= \dfrac{2}{N^{6}} \sum_{j<k}^{N-1} \left[ (j+3)(j+2)(k+6)(k+5  \right] \nonumber
\end{align}
and this  can be written as  
\begin{align}
Q_{N}^{(IV)} - Q_{N}^{(V)} - Q_{N}^{(VI)} + Q_{N}^{(VII)}  &= \dfrac{2}{N^{6}} \sum_{k=1}^{N_{\alpha}-1} (k+6)(k+5) \sum_{j=0}^{k-1} (j+3)(j+2)\nonumber \\
&\sim \dfrac{2}{N^{6}} \sum_{k=1}^{N_{\alpha}-1} (k+6)(k+5)  k^ {3} + o\left( \frac{1}{N} \right) \nonumber \\
&  \sim \frac{2}{18} \frac{(N_{\alpha})^ {6}}{N ^ {6} } +  o\left( \frac{1}{N} \right) \nonumber \\ 
&\xrightarrow[N \to \infty]{} \dfrac{\alpha^6}{9}. \label{lim2}
\end{align}

Finally, considering the results obtained in  \eqref{lim1} and \eqref{lim2}, we can conclude.
\qed

\subsection{Joint densities under Renewal Sampling}
\begin{table}[h!]
\centering
\begin{tabular}{|c|c|c|}
\hline
Joint distribution & Probability Density Function & Support \\\hline
$f_{\tau_{i} , t_{i+1}}({a, b})$ & $\dfrac{N^{i+1}}{\Gamma(i)} a^{i-1} e^{-N (a + b)}$ & $0 \leq a < \infty$ \\
                                &                                                                                                          & $0 \leq b < \infty$ \\ \hline
$f_{\tau_{i} , \tau_{i+1}}({a, b}) $ & $\dfrac{N^{i+1}}{\Gamma(i)} a^{i-1}  e^{-N b}$ & $0 \leq a \leq b$  \\
                                  &                                                                                          & $0 \leq b < \infty$ \\ \hline
$f_{\tau_{i-1} , t_{i}, t_{i+1}}({a, b, c})$ & $\dfrac{N^{i+1}}{\Gamma(i-1)} a^{i-2} e^{-N (a + b + c)}$ & $0 \leq a < \infty $ \\
                                                &                                                                                                                                     & $0 \leq b < \infty$ \\ 
                                                &                                                                                                                                     &  $0 \leq c < \infty$ \\ \hline
$f_{\tau_{i-1} , \tau_{i} , \tau_{i+1}}({a, b, c})$ & $\dfrac{N^{i+1}}{\Gamma(i-1)} a^{i-2} e^{-N c}$ & $0 \leq a \leq b$ \\
                                              &                                                                                                   & $0 \leq b \leq c$ \\
                                              &                                                                                                   & $0 \leq c < \infty$ \\ \hline
$f_{\tau_{j} , \tau_{j+1} , \tau_{i} , \tau_{i+1}}({a, b, c,d})$ & $\dfrac{N^{i+1}}{\Gamma(j) \Gamma(i-j-1)} a^{j-1} \left(c - b  \right)^{i-j-2} e^{-N d}$ & $0 \leq a \leq b$\\ 
                                                           &                                                                                      & $0 \leq b \leq c$ \\
                                                           &                                                                                      & $0 \leq c \leq d$ \\
                                                           &                                                                                      &  $0 \leq d < \infty$ \\ \hline
\end{tabular}
\caption{Densities under Renewal Process} \label{densities}
\end{table}

\noindent {\bf Acknowledgements:}
This  research was partially supported by Project REDES 150038, MATHAMSUD 19-MATH-06, Math AmSud 18-MATH-07 SaSMoTiDep, CONICYT - MATHAMSUD FANTASTIC  20-MATH-05. T. Roa was partially supported by Beca CONICYT-PFCHA/Doctorado Nacional/2018-21180298, S. Torres was partially supported by FONDECYT 1171335 and C. Tudor was partially supported by MEC PAI80160046.
\par

\end{document}